\documentclass[11pt]{article}
\usepackage{amssymb,amsfonts,amsmath}
\usepackage{epsfig}
\usepackage{graphicx}
\usepackage{color}
\usepackage{anysize}
\usepackage{hyperref}

\marginsize{2.3cm}{2.3cm}{1cm}{2cm}

\begin{document}
\baselineskip=14pt

\makeatletter\@addtoreset{equation}{section}
\makeatother\def\theequation{\thesection.\arabic{equation}}

\newtheorem{defin}{Definition}[section]
\newtheorem{Prop}{Proposition}
\newtheorem{teo}{Theorem}[section]
\newtheorem{ml}{Main Lemma}
\newtheorem{con}{Conjecture}
\newtheorem{cond}{Condition}
\newtheorem{conj}{Conjecture}
\newtheorem{prop}[teo]{Proposition}
\newtheorem{lem}{Lemma}[section]
\newtheorem{rmk}[teo]{Remark}
\newtheorem{cor}{Corollary}[section]

\newcommand{\be}{\begin{equation}}
\newcommand{\ee}{\end{equation}}
\newcommand{\ben}{\begin{eqnarray}}
\newcommand{\benn}{\begin{eqnarray*}}
\newcommand{\een}{\end{eqnarray}}
\newcommand{\eenn}{\end{eqnarray*}}
\newcommand{\bp}{\begin{prop}}
\newcommand{\ep}{\end{prop}}
\newcommand{\bt}{\begin{teo}}
\newcommand{\et}{\end{teo}}
\newcommand{\bcor}{\begin{cor}}
\newcommand{\ecor}{\end{cor}}
\newcommand{\bcon}{\begin{con}}
\newcommand{\econ}{\end{con}}
\newcommand{\bcond}{\begin{cond}}
\newcommand{\econd}{\end{cond}}
\newcommand{\br}{\begin{rmk}}
\newcommand{\er}{\end{rmk}}
\newcommand{\bl}{\begin{lem}}
\newcommand{\el}{\end{lem}}
\newcommand{\bit}{\begin{itemize}}
\newcommand{\eit}{\end{itemize}}
\newcommand{\bd}{\begin{defin}}
\newcommand{\ed}{\end{defin}}
\newcommand{\bpr}{\begin{proof}}
\newcommand{\epr}{\end{proof}}

\newenvironment{proof}{\noindent {\em Proof}.\,\,}{\hspace*{\fill}$\halmos$\medskip}

\newcommand{\halmos}{\rule{1ex}{1.4ex}}
\def \qed {{\hspace*{\fill}$\halmos$\medskip}}

\newcommand{\Z}{{\mathbb Z}}
\newcommand{\R}{{\mathbb R}}
\newcommand{\E}{{\mathbb E}}
\newcommand{\C}{{\mathbb C}}
\renewcommand{\P}{{\mathbb P}}
\newcommand{\N}{{\mathbb N}}

\newcommand{\Bi}{{\cal B}}
\newcommand{\Si}{{\cal S}}
\newcommand{\Ti}{{\cal T}}
\newcommand{\Wi}{{\cal W}}
\newcommand{\Yi}{{\cal Y}}
\newcommand{\Hi}{{\cal H}}
\newcommand{\Fi}{{\cal F}}
\newcommand{\Zi}{{\cal Z}}

\newcommand{\eps}{\epsilon}

\newcommand{\nn}{\nonumber}

\newcommand{\pa}{\partial}
\newcommand{\ffrac}[2]{{\textstyle\frac{{#1}}{{#2}}}}
\newcommand{\dif}[1]{\ffrac{\partial}{\partial{#1}}}
\newcommand{\diff}[1]{\ffrac{\partial^2}{{\partial{#1}}^2}}
\newcommand{\difif}[2]{\ffrac{\partial^2}{\partial{#1}\partial{#2}}}

\title{Disorder relevance for the random walk pinning model in dimension $3$}
\author{Matthias Birkner$^{\,1}$, Rongfeng Sun$^{\,2}$}
\date{Apr 20, 2010}
\maketitle

\footnotetext[1]{Institut f\"ur Mathematik, Johannes-Gutenberg-Universit\"at Mainz, Staudingerweg 9, 55099 Mainz, Germany. Email:
birkner@mathematik.uni-mainz.de}

\footnotetext[2]{Department of Mathematics, National University of Singapore, 10 Lower Kent Ridge Road, 119076 Singapore. Email:
matsr@nus.edu.sg}


\begin{abstract}
We study the continuous time version of the {\it random walk pinning model}, where conditioned on a continuous time
random walk $(Y_s)_{s\geq 0}$ on $\Z^d$ with jump rate $\rho>0$, which plays the role of disorder, the law up to time
$t$ of a second independent random walk $(X_s)_{0\leq s\leq t}$ with jump rate $1$ is Gibbs transformed with weight
$e^{\beta L_t(X,Y)}$, where $L_t(X,Y)$ is the collision local time between $X$ and $Y$ up to time $t$. As the inverse
temperature $\beta$ varies, the model undergoes a localization-delocalization transition at some critical $\beta_c\geq 0$.
A natural question is whether or not there is disorder relevance, namely whether or not $\beta_c$ differs from the critical
point $\beta^{\rm ann}_c$ for the annealed model. In \cite{BS09}, it was shown that there is disorder irrelevance
in dimensions $d=1$ and $2$, and disorder relevance in $d\geq 4$. For $d\geq 5$, disorder relevance was first
proved in \cite{BGdH08}. In this paper, we prove that if $X$ and $Y$ have the same jump probability kernel, which
is irreducible and symmetric with finite second moments, then there is also disorder relevance in the critical
dimension $d=3$, and $\beta_c-\beta^{\rm ann}_c$ is at least of the order $e^{-C(\zeta)/\rho^\zeta}$, $C(\zeta)>0$, for any $\zeta>2$.
Our proof employs coarse graining and fractional moment techniques, which have recently been applied by
Lacoin \cite{L09} to the directed polymer model in random environment, and by Giacomin, Lacoin and Toninelli \cite{GLT09}
to establish disorder relevance for the random pinning model in the critical dimension. Along the way, we also prove a
continuous time version of Doney's local limit theorem \cite{D97} for renewal processes with infinite mean.


\medskip


\end{abstract}

\section{Model and result}

Let us recall the continuous time {\it random walk pinned to random walk model} studied in \cite{BS09}, which we will abbreviate from now on as the
{\it random walk pinning model} (RWPM). Let $X$ and $Y$ be two continuous time random walks on $\Z^3$ starting from the origin, such that $X$ and $Y$ have respectively jump rates $1$ and $\rho\geq 0$, and identical irreducible symmetric jump probability kernels on $\Z^3$ with finite second moments. Let $\mu_t$ denote the law of $(X_s)_{0\leq s\leq t}$. Then given $\beta\in\R$ and conditioned on $(Y_s)_{s\geq 0}$, which we interpret as a random environment or disorder, we define a Gibbs transform $\mu^\beta_{t,Y}$ of the path measure $\mu_t$ via Radon-Nikodym derivative
\be\label{mucontinuous}
\frac{{\rm d} \mu^\beta_{t,Y}}{{\rm d\mu_t}}(X) = \frac{e^{\beta L_t(X,Y)}}{Z^\beta_{t,Y}},
\ee
where $L_t(X,Y)=\int_0^t 1_{\{X_s=Y_s\}}{\rm d}s$, and
\be\label{Zcts}
Z^\beta_{t,Y} = \E^X_0\big[e^{\beta L_t(X,Y)}\big]
\ee
is the {\it quenched partition function} with $\E^X_x[\cdot]$ denoting expectation w.r.t.\ $X$ starting from $x\in\Z^3$. We can interpret $X$ as a polymer which is
attracted to a random defect line $Y$. A more commonly studied model is to consider a constant defect line $Y\equiv 0$, but with random strength of interaction
between $X$ and $Y$ at different time points. This is known as the {\em random pinning model} (RPM), the discrete time analogue of which was the subject of many recent papers
(see e.g.~\cite{DGLT09, GLT08, GLT09}), as well as a book~\cite{G07}.

A common variant of the Gibbs measure $\mu^\beta_{t,Y}$ is to introduce pinning of path at the end point $t$, i.e., we define the Gibbs measure $\mu^{\beta,\rm pin}_{t,Y}$ with
\be\label{mupincts}
\frac{{\rm d}\mu^{\beta, \rm pin}_{t,Y}}{{\rm d}\mu_t}(X) = 1_{\{X_t=Y_t\}} \frac{e^{\beta L_t(X,Y)}}{Z^{\beta, \rm pin}_{t,Y}}
\ee
with $Z^{\beta, \rm pin}_{t,Y} = \E^X_0\big[e^{\beta L_t(X,Y)}1_{\{X_t=Y_t\}}\big]$. It was shown in \cite{BS09} that, almost surely w.r.t.\ $Y$, the limit
\be\label{Fbetarho}
F(\beta, \rho) = \lim_{t\to\infty} \frac{1}{t} \log Z^\beta_{t,Y} = \lim_{t\to\infty} \frac{1}{t} \log Z^{\beta,\rm pin}_{t,Y}
\ee
exists and is independent of the disorder $Y$, which we call the {\it quenched free energy} of the model. There exists a critical inverse temperature $\beta_c=\beta_c(\rho)$, such that $F(\beta,\rho)>0$ if $\beta>\beta_c$ and $F(\beta)=0$ if $\beta<\beta_c$. The supercritical region $\beta \in (\beta_c,\infty)$
is the localized phase where given $Y$, and with respect to either $\mu^{\beta}_{t,Y}$ or $\mu^{\beta,\rm pin}_{t,Y}$, the contact fraction $L_t(X,Y)/t$ between $X$
and $Y$ up to time $t$ typically remains positive as $t\to\infty$, so that the walk $X$ is pinned to $Y$. In fact, by the convexity of $\log Z^\beta_{t,Y}$ in $\beta$ and (\ref{Fbetarho}), it is not hard to see that almost surely,
\begin{eqnarray*}
\liminf_{t\to\infty} \mu^{\beta}_{t,Y}(t^{-1} L_t(X,Y)) &=& \liminf_{t\to\infty} \frac{\partial (t^{-1}\log Z^\beta_{t,Y})}{\partial \beta} \geq \frac{\partial F(\beta,\rho)}{\partial_-\beta}, \\
\limsup_{t\to\infty} \mu^{\beta}_{t,Y}(t^{-1} L_t(X,Y)) &=& \limsup_{t\to\infty} \frac{\partial (t^{-1}\log Z^\beta_{t,Y})}{\partial \beta} \leq \frac{\partial F(\beta,\rho)}{\partial_+\beta},
\end{eqnarray*}
where $\frac{\partial}{\partial_-\beta}$ and $\frac{\partial}{\partial_+\beta}$ denote respectively the left and right derivative w.r.t.\ $\beta$. The convexity
of $F(\beta,\rho)$ in $\beta$ implies that $\frac{\partial F(\beta,\rho)}{\partial_-\beta}>0$ for all $\beta>\beta_c$. In contrast, the subcritical region $\beta \in (-\infty, \beta_c)$ is the de-localized phase, where $\partial F(\beta,\rho)/\partial \beta=0$ and the contact fraction $L_t(X,Y)/t$ is typically of order $o(1)$ as $t\to\infty$, so that $X$ becomes delocalized from $Y$.

An important tool in the study of models with disorder is to compare the quenched free energy with the annealed free energy, which is defined by
\be\label{Fanncts}
F_{\rm ann}(\beta,\rho) := \lim_{t\to\infty} \frac{1}{t} \log Z^{\beta}_{t,\rm ann} = \lim_{t\to\infty}\frac{1}{t} \log Z^{\beta, \rm pin}_{t,\rm ann},
\ee
where
$$
Z^{\beta}_{t,\rm ann} = \E^{X,Y}_{0,0}[e^{\beta L_t(X,Y)}] \qquad \mbox{and} \qquad Z^{\beta, \rm pin}_{t,\rm ann} = \E^{X,Y}_{0,0}[e^{\beta L_t(X,Y)}1_{\{X_t=Y_t\}}]
$$
are the free (resp.\ constrained) versions of the annealed partition function for the RWPM. Since $X-Y$ is also a random walk, we see that
$Z^{\beta}_{t,\rm ann}$ and $Z^{\beta, \rm pin}_{t,\rm ann}$ are the partition functions of a RWPM where the random walk $X-Y$ is attracted
to the constant defect line $0$. This defines the annealed model. In particular, there also exists a critical
point $\beta^{\rm ann}_c=\beta^{\rm ann}_c(\rho)$ such that $F_{\rm ann}(\beta,\rho)>0$ when $\beta>\beta^{\rm ann}_c$ and $F_{\rm ann}(\beta,\rho)=0$ when $\beta<\beta^{\rm ann}_c$. It is easy to show that $\beta_c^{\rm ann} = (1+\rho)/G$, where $G$ is the Green function of $X$, see end of Sec.~\ref{S:rep},
while no explicit expression for $\beta^c$ is known. By Jensen's inequality, it is easily seen that $F(\beta,\rho) \leq F_{\rm ann}(\beta,\rho)$, and hence $\beta_c\geq \beta^{\rm ann}_c$. A fundamental question is then to determine whether the disorder is sufficient to shift the critical point of the model so that $\beta_c > \beta^{\rm ann}_c$, which is called {\it disorder relevance}. If $\beta_c=\beta^{\rm ann}_c$, then we say there is {\it disorder irrelevance}, and it is generally believed that the quenched model's
behavior in this case is similar to that of the annealed model. It turns out that disorder relevance/irrelevance has an interesting dependence on the spatial dimension $d$.

In \cite{BS09}, it was shown that if $X$ and $Y$ are continuous time simple random walks, then the RWPM is disorder irrelevant in $d=1$ and $2$, and disorder relevant in $d\geq 4$. Furthermore, it was shown that in $d\geq 5$, there exists $a>0$ such that $\beta_c-\beta^{\rm ann}_c>a\rho$ for all
$\rho\in [0,1]$; while in $d=4$, for any $\delta>0$, there exists $a_\delta>0$ such that $\beta_c-\beta_c^{\rm ann}\geq a_\delta \rho^{1+\delta}$ for all $\rho\in [0,1]$.
It is easy to check that the proof of these results in \cite{BS09} apply equally well to continuous time random walks $X$ and $Y$ with the same irreducible symmetric
jump probability kernel with finite second moments. In this paper, we resolve the marginal case $d=3$ and show that there is disorder relevance.

\bt{\bf [Annealed vs quenched critical points]}\label{T:cptcts}\\
Let $X$ and $Y$ be two continuous time random walks with respective jump rates $1$ and $\rho>0$ and identical irreducible symmetric jump probability kernel $q(\cdot)$ on
$\Z^3$ with finite second moments. Assume $X_0=Y_0=0$. Then for the associated RWPM, $\beta_c(\rho)>\beta^{\rm ann}_c(\rho)$ for all $\rho>0$, and for any $\zeta>2$, there exists $C(\zeta)>0$ such that for all $\rho \in (0,1]$,
\be\label{eq:cptcts}
\beta_c(\rho)-\beta^{\rm ann}_c(\rho)\geq e^{-C(\zeta)\rho^{-\zeta}}.
\ee
\et
{\bf Remark.} It is intriguing that our lower bound for the critical point shift is of the same form as for the RPM in the marginal case, where
a lower bound of $e^{-C(\zeta)\beta^{-\zeta}}$ was obtained in \cite{GLT09} for any $\zeta>2$, and $\zeta=2$ is known to provide an upper bound. For the RWPM, there has
been no heuristics or results so far on the upper bound.
\medskip

Let
\be\label{beta*cts}
\beta^*_{c}(\rho) = \sup\Big\{\beta \in\R : \sup_{t>0} Z^{\beta}_{t,Y} <\infty\ \mbox{a.s.  w.r.t. } Y\Big\}.
\ee
Note that $\beta^*_c(\rho) \leq \beta_c(\rho)$. We will in fact prove the following stronger version of Theorem \ref{T:cptcts}.

\bt{\bf [Non-coincidence of critical points strengthened]}\label{T:cptctsst}\\
Assuming the same conditions as in Theorem \ref{T:cptcts}, then the conclusions therein also hold with $\beta_c(\rho)$ replaced by $\beta^*_c(\rho)$.
\et
{\bf Remark.} The question whether $\beta^*_c=\beta_c$ or $\beta^*_c<\beta_c$ remains open, and so is the analogous question for the RPM. Note that when $Z^\beta_{t,Y}$
is uniformly bounded in $t>0$, the distribution of $L_t(X,Y)$ under the measure $\mu^\beta_{t,Y}$ remains tight as $t\to\infty$, while $Z^\beta_{t,Y}\to\infty$ if and only if $L_t(X,Y)$ under $\mu^\beta_{t,Y}$ tends to $\infty$ in probability. If $\beta^*_c<\beta_c$, then there exists a phase in the delocalized regime where
$L_t(X,Y)$ under $\mu^\beta_{t,Y}$ tends to $\infty$ at a rate that is $o(t)$, which would be very surprising.
\bigskip

Theorem \ref{T:cptctsst} confirms a conjecture of Greven and den Hollander \cite[Conj.~1.8]{GdH07} that in $d=3$, the Parabolic Anderson Model (PAM) with Brownian noise could admit an equilibrium measure with an infinite second moment. We refer to \cite[Sec.~1.4]{BS09} for a more detailed discussion on the connection between the RWPM and the PAM, as well as the connection of the discrete time RPMs and RWPMs with the directed polymer model in random environment.

Our proof of Theorem~\ref{T:cptctsst} will follow the general approach developed by Giacomin, Lacoin, and Toninelli in~\cite{GLT08, GLT09} for proving the marginal relevance of disorder for the random pinning model (RPM), as well as by Lacoin in~\cite{L09} for the study of the directed polymer model in random environment. The basic ingredients are change of measure arguments for bounding fractional moments of the partition function $Z^\beta_{t,Y}$, coupled with a coarse grain splitting of $Z^\beta_{t,Y}$. These techniques have proven to be remarkably powerful, and they apply to a wide
range of models: in particular, to weighted renewal processes in random environments, including the random pinning, the random walk pinning, and the copolymer models
(see \cite[Sec.~1.4]{BS09} for a more detailed discussion), as well as to weighted random walks in random environments, including the directed polymer model \cite{L09}
and random walk in random environments~\cite{YZ09}. We will recall in detail the fractional moment techniques and the coarse graining procedure and formulate them for the RWPM, which will constitute the model independent part of our analysis. A key element of the fractional moment
technique involves a change of measure, and more generally, the choice of a suitable test function. This is the model dependent part of the analysis, which in general
is far from trivial since disorder relevance in the critical dimension is a rather subtle effect. The bulk of this paper is thus dedicated to the choice of a suitable test function for the RWPM and its analysis. Compared to the RPM and the directed polymer model, new complications arise due to the different nature of the disorder
of the RWPM.
\medskip

We also include here a result on the monotonicity of $\beta_c(\rho)-\beta^{\rm ann}_c(\rho)$, resp.\ $\beta^*_c(\rho)-\beta^{\rm ann}_c(\rho)$, in $\rho$, which was pointed out to us by the referee along with an elegant proof.
\bt{\bf [Monotonicity of critical point shift]}\label{T:monshif} \\
Assuming the same conditions as in Theorem \ref{T:cptcts} for a RWPM in $\Z^d$ with $d\geq 3$, we have
\be\label{monoton1}
\left.
\begin{aligned}
\frac{\beta_c(\rho')}{1+\rho'}\geq \frac{\beta_c(\rho)}{1+\rho}, \\
\frac{\beta^*_c(\rho')}{1+\rho'}\geq \frac{\beta^*_c(\rho)}{1+\rho},
\end{aligned}
\right.
\qquad \mbox{for all } \rho'>\rho\geq 0.
\ee
In particular,
\be\label{monoton2}
\left.
\begin{aligned}
\beta_c(\rho')-\beta^{\rm ann}_c(\rho') > \beta_c(\rho)-\beta^{\rm ann}_c(\rho), \\
\beta^*_c(\rho')-\beta^{\rm ann}_c(\rho') > \beta^*_c(\rho)-\beta^{\rm ann}_c(\rho),
\end{aligned}
\right.
\qquad \mbox{for all } \rho'>\rho\geq 0.
\ee
\et
We defer its proof to Appendix~\ref{S:monshif}. We remark that proving the strict inequalities in (\ref{monoton2}) requires Theorem~\ref{T:cptctsst} and its analogue in dimensions $d\geq 4$.
\medskip

{\bf Outline.}
The rest of the paper is organized as follows. In Section \ref{S:rep}, we recall from \cite{BGdH08} and \cite{BS09} a representation of $Z^{\beta}_{t,Y}$ as
the partition function of a weighted renewal process in random environment. In Section \ref{S:fraction}, we recall the coarse graining procedure and fractional moment techniques developed in \cite{GLT08}, \cite{L09} and \cite{GLT09}. To prove Theorem \ref{T:cptctsst}, we apply the coarse graining procedure to $Z^{\beta}_{t,Y}$ instead of the constrained partition function $Z^{\beta, \rm pin}_{t,Y}$ as done in \cite{GLT08, GLT09}. The proof of disorder relevance is then reduced in Section \ref{S:fraction} to two key propositions: Prop.~\ref{P:cgZab}, which is model dependent and needs to be proved for any new weighted renewal process in random environment one is interested in, and Prop.~\ref{P:coarsegrain}, which is model independent.  Compared to analogues of Prop.~\ref{P:cgZab} formulated previously for the RPM (see \cite[Lemma 3.1]{GLT09}), our weaker formulation here (more precisely its reduction to Prop.~\ref{P:cgZfree} in Sec.~\ref{S:cgZab}) allows a more direct comparison with renewal processes without boundary constraints, which conceptually simplifies subsequent analysis. In Section \ref{S:HLY}, we identify a crucial test function $H_L(Y)$ for the disorder $Y$ and state some essential properties. Assuming these properties, we then prove in Section \ref{S:cgZab} the key Prop.~\ref{P:cgZab}, which is further reduced to a model dependent Prop.~\ref{P:cgZfree} by extracting some model independent renewal calculations. In Section \ref{S:coarsegrain}, we deduce Prop.~\ref{P:coarsegrain} from Prop.~\ref{P:cgZab}, which is again model independent. The properties of the test function $H_L$ are then established in Sections \ref{S:HLYbound}--\ref{S:Hsigmabd}. In Appendices \ref{S:renewal} and \ref{S:RW}, we prove some renewal and random walk estimates which we need for our proof. In particular, we prove in Lemma \ref{L:renewalprob} a continuous time version of Doney's local limit theorem \cite[Thm.~3]{D97} for renewal processes with infinite mean.
Finally, in Appendix \ref{S:monshif}, we include a proof of Theorem~\ref{T:monshif} shown to us by the referee.
\bigskip

{\bf Note.} During the preparation of this manuscript, we became aware of a preprint by Q.~Berger and F.L.~Toninelli~\cite{BT09}, in which
they proved disorder relevance for the {\em discrete time} RWPM in dimension $3$ under the assumption that the random walk increment is
symmetric with sub-Gaussian tails. An inspection shows that the main difference between our two approaches lies in the choice of the test function
$H_L(\cdot)$ in (\ref{HLY}), which results in completely different model dependent analysis as well as different assumptions on the model. In principle, both
approaches should be applicable to both discrete and continuous time models. Most results in this paper carry over directly to the discrete time case. The only exception
is Lemma \ref{L:retprocomp}, for which we do not have a proof for its discrete time analogue. Lemma \ref{L:retprocomp} is used to prove Lemma \ref{L:Hsigmabd}~(\ref{Hsigmext})--(\ref{Hsigmext2}). In light of \cite{BT09}, we will not pursue this further in this paper.
\bigskip

{\bf Notation:} Throughout the rest of this paper, unless stated otherwise, we will use $C$, $C_1$ and $C_2$ to denote generic constants whose precise values may
change from line to line. However, their values all depend only on the jump rate $\rho$ and the jump probability kernel $q(\cdot)$, and are uniform in $\rho\in (0,1]$.

\section{Representation as a weighted renewal process in random environment}\label{S:rep}

First we recall from \cite{BGdH08} and \cite{BS09} a representation of $Z^\beta_{t,Y}$ as the partition function of a weighted renewal process in random environment.
Let $p_s(\cdot)=p^X_s(\cdot)$ denote the transition probability kernel of $X$ at time $s$. Then $Y$ and $X-Y$ have respective transition kernels
$p^Y_s(\cdot):=p_{\rho s}(\cdot)$ and $p^{X-Y}_s(\cdot):=p_{(1+\rho)s}(\cdot)$. Let
$$
G=\int_0^\infty p_s(0){\rm d}s, \quad G^{X-Y}=\int_0^\infty p_{(1+\rho)s}(0){\rm d}s=\frac{G}{1+\rho}, \quad K(t) = \frac{p^{X-Y}_t(0)}{G^{X-Y}} = \frac{(1+\rho)p_{(1+\rho)t}(0)}{G},
$$
where $K(t)dt$ is to be interpreted as the renewal time distribution of a recurrent renewal process $\sigma=\{\sigma_0=0<\sigma_1<\cdots\} \subset [0,\infty)$. Let $z=\beta G^{X-Y}=\beta G/(1+\rho)$ and $\Zi^z_{t,Y}:= Z^\beta_{t,Y}$. Then
\begin{eqnarray}
\Zi^{z}_{t,Y} &=& \E^X_0\big[e^{\beta L_t(X,Y)} \big] =
\E^X_0\left[1+\sum_{m=1}^\infty \frac{\beta^m}{m!}\Big(\int_0^t 1_{\{X_s=Y_s\}}ds\Big)^m\right]  \nonumber\\
&=&
\E^X_0\Big[1+\sum_{m=1}^\infty \beta^m \idotsint\limits_{0<\sigma_1\cdots<\sigma_m<t}
1_{\{X_{\sigma_1}=Y_{\sigma_1},\cdots, X_{\sigma_m}=Y_{\sigma_m}\}}d\sigma_1\cdots d\sigma_m\Big]  \nonumber \\
&=&
1 + \sum_{m=1}^\infty \beta^m \idotsint\limits_{0<\sigma_1\cdots<\sigma_m<t}
p_{\sigma_1}(Y_{\sigma_1})p_{\sigma_2-\sigma_1}(Y_{\sigma_2}-Y_{\sigma_1})\cdots p_{\sigma_m-\sigma_{m-1}}(Y_{\sigma_m}-Y_{\sigma_{m-1}})d\sigma_1\cdots d\sigma_m \nonumber \\
&=&
1 + \sum_{m=1}^\infty z^m \!\!\!\!\!
\idotsint\limits_{ \atop \sigma_0=0<\sigma_1\cdots<\sigma_m<t}
\prod\limits_{i=1}^{m} K(\sigma_i-\sigma_{i-1}) W(\sigma_i-\sigma_{i-1}, Y_{\sigma_i}-Y_{\sigma_{i-1}})\ d\sigma_1\cdots
d\sigma_m, \label{ctszfree}
\end{eqnarray}
where
\be\label{W}
W(\sigma_i-\sigma_{i-1}, Y_{\sigma_i}-Y_{\sigma_{i-1}})=
\frac{p^X_{\sigma_i-\sigma_{i-1}}(Y_{\sigma_i}-Y_{\sigma_{i-1}})}{p^{X-Y}_{\sigma_i-\sigma_{i-1}}(0)}.
\ee
We can thus interpret $\Zi^z_{t,Y}$ as the partition function of a weighted renewal process $\sigma$ in the random environment $Y$, where the renewal
time distribution is given by $K(\cdot)$, and the $i$-th renewal return incurs a weight factor of $zW(\sigma_i-\sigma_{i-1}, Y_{\sigma_i}-Y_{\sigma_{i-1}})$.

Similarly, for any $0\leq U \leq V$, we can define $\Zi^{z, {\rm pin}}_{[U,V],Y}:=1$ when $U=V$, and otherwise
\begin{eqnarray}
\Zi^{z, {\rm pin}}_{[U,V],Y} := \!\! \sum_{m=1}^\infty \!\!
\idotsint\limits_{ \atop \sigma_0=U<\sigma_1\cdots<\sigma_m=V} \!\!\!\!\!\!\!
z^m  \prod\limits_{i=1}^{m} K(\sigma_i-\sigma_{i-1}) W(\sigma_i-\sigma_{i-1}, Y_{\sigma_i}-Y_{\sigma_{i-1}})\ d\sigma_1\cdots
d\sigma_{m-1}, \label{ctsz}
\end{eqnarray}
where the term for $m=1$ is defined to be $zK(V-U)W(V-U,Y_V-Y_U)$. Note that $\Zi^{z,\rm pin}_{[0,t],Y} = \beta Z^{\beta, \rm pin}_{t,Y}$, which we will
simply denote by $\Zi^{z,\rm pin}_{t,Y}$.

Since $K$ is the renewal time distribution of a recurrent renewal process $\sigma$ on $[0,\infty)$, and note that
$\E^Y_0[W(v-u, Y_v-Y_u)]=1$ for any $u<v$, the critical point $z^{\rm ann}_c$ of the annealed model with partition function $\E^Y_0[\Zi^z_{t,Y}]$
is exactly $1$. By the mapping $z= \beta G^{X-Y}$, we deduce that $\beta^{\rm ann}_c=1/G^{X-Y}=(1+\rho)/G$.
The mapping to a weighted renewal process in random environment casts the RWPM in the same framework as the RPM, which paves the way for the application
of general approaches developed in \cite{GLT08, GLT09}.

\section{Fractional moment techniques and coarse graining}\label{S:fraction}

We now recall the fractional moment techniques and the coarse graining procedure, which were developed in a series of papers for the RPM that culminated
in \cite{GLT08, GLT09}, where marginal relevance of disorder was established, as well as in \cite{L09} where the same techniques
were applied to the directed polymer model in random environment.

By (\ref{Zcts}), $Z^\beta_{t,Y}=\Zi^z_{t,Y}$ is monotonically increasing in $t$ for every realization of $Y$. Therefore, to prove Theorem \ref{T:cptctsst},
it suffices to show that for some $\gamma \in (0,1)$, and some $z>1$ depending suitably on $\rho$, we have
\be\label{Zgammabd}
\sup_{t>0} \E^Y_0[ (\Zi^{z}_{t,Y})^\gamma]<\infty.
\ee
We will choose below a coarse graining scale $L$ and show that for each $\rho>0$, (\ref{Zgammabd}) holds for all $z\in (1, 1+1/L]$ if $L$ is
sufficiently large, which implies $z_c\geq 1+1/L$ and $\beta_c^*-\beta^{\rm ann}_c\geq (1+\rho)/GL$. For $\rho \in (0,1]$, we will let
$L=e^{B_1/\rho^\zeta}$ for any fixed $\zeta>2$, and prove that (\ref{Zgammabd}) holds for all $z\in (1, 1+1/L]$ uniformly in $\rho\in (0,1]$ if $B_1$ is
large enough. This would then imply the lower bound on $\beta_c^*-\beta_c^{\rm ann}$ in Theorem \ref{T:cptctsst}.

Note that by using more refined large deviation estimates for the renewal process
with waiting time density $K(\cdot)$, it seems possible to extend (\ref{Zgammabd}) to
$z \in (1, 1+1/L^\eta]$ with a suitable $\eta \in (1/2, 1)$. By the relation between $L$
and $\rho$ in our coarse-graining scheme, this would only affect the (unspecified)
constant $C(\zeta)$ in (\ref{eq:cptcts}), not the exponent $\zeta$ itself.

To bound the fractional moment $\E^Y_0[ (\Zi^{z}_{t,Y})^\gamma]$, we apply the inequality
\be\label{fracmomineq}
(\sum_{i=1}^n a_i)^\gamma \leq \sum_{i=1}^n a_i^\gamma \quad \quad \mbox{for}\quad a_i\geq 0,\ 1\leq i\leq n, \quad \mbox{and} \quad \gamma \in (0,1).
\ee
This seemingly trivial inequality turns out to be exceptionally powerful in bounding fractional moments.
However, the success of such a bound depends crucially on how $\Zi^z_{t,Y}$ is split into a sum of terms. This is where coarse graining comes in,
which was used in \cite{GLT08, L09, GLT09}. We remark that in the earlier paper \cite{DGLT09} on the RPM, and later in the analysis~\cite{BS09}
of the RWPM in $d\geq 4$ , $\Zi^z_{t,Y}$ is partitioned according to the values of the pair of consecutive renewal times $\sigma_i<\sigma_{i+1}$
which straddle a fixed time $L>0$. The coarse graining procedure we recall below uses a more refined partition of $\Zi^z_{t,Y}$.

Fix a large constant $L>0$, which will be the coarse graining scale. Assume that $t=mL$ for some $m\in\N$. Then we partition $(0,t]$
into $m$ blocks $\Lambda_1,\cdots, \Lambda_m$ with $\Lambda_i:=((i-1)L,iL]$. The coarse graining procedure simply groups terms in (\ref{ctszfree}) according to which blocks $\Lambda_i$
does the renewal configuration $\sigma:=\{\sigma_0=0<\sigma_1<\cdots\}$ intersect. More precisely, the set of blocks in $\{\Lambda_i\}_{1\leq i\leq m}$ which $\sigma$
intersects can be represented by a set $I\subset \{1,\cdots, m\}$. Then we can decompose $\Zi^z_{t,Y}$ in (\ref{ctszfree}) as
$$
\Zi^z_{t,Y} = \sum_{I\subset \{1,\cdots, m\}} \Zi^{z,I}_{t,Y},
$$
where $\Zi^{z,\emptyset}_{t,Y}:=1$, and for each $I=\{1\leq i_1 < i_2<\cdots < i_k\leq m\} \neq \emptyset$,
\be \label{ctsZzINY}
\Zi^{z,I}_{t,Y} =  \!\!\!\!\! \int\limits_{a_1<b_1 \atop a_1,b_1\in \Lambda_{i_1}}\!\!\!\!\!\cdots\!\!\!\!\! \int\limits_{a_k< b_k\atop a_k,b_k\in \Lambda_{i_k}} \!\!\!\!\!
\prod_{j=1}^k K(a_j-b_{j-1}) z W(a_j-b_{j-1}, Y_{a_j}-Y_{b_{j-1}}) \Zi^{z,\rm pin}_{[a_j, b_j],Y} \prod_{j=1}^k {\rm d}a_j\, {\rm d}b_j,
\ee
where $b_0:=0$. By (\ref{fracmomineq}), for any $\gamma\in (0,1)$, we have
\be\label{Zgsplit}
\E^Y_0\big[(\Zi^z_{t,Y})^\gamma\big] \leq \sum_{I\subset \{1,\cdots, m\}} \E^Y_0\big[(\Zi^{z,I}_{t,Y})^\gamma\big].
\ee
We will prove (\ref{Zgammabd}) by comparing $\E^Y_0\big[(\Zi^{z,I}_{t,Y})^\gamma\big]$ with the probability that a subcritical renewal
process on $\N\cup\{0\}$ intersects $\{1,\cdots, m\}$ exactly at $I$.

To bound $\E^Y_0\big[(\Zi^{z,I}_{t,Y})^\gamma\big]$, one introduces a change of measure. Let $f_I(Y)$ be a non-negative function of the disorder $Y$. By
H\"older's inequality,
\be\label{Holder}
\E^Y_0\big[\big(\Zi^{z, I}_{t,Y}\big)^\gamma\big] = \E^Y_0\big[f_I(Y)^\gamma f_I(Y)^{-\gamma} \big(\Zi^{z, I}_{t,Y}\big)^\gamma\big]
\leq \E^Y_0\big[f_I(Y)^{-\frac{\gamma}{1-\gamma}}\big]^{1-\gamma} \E^Y_0\big[f_I(Y) \Zi^{z,I}_{t,Y}\big]^\gamma.
\ee
To decouple different blocks $\Lambda_i$, we will let $f_I(Y) = \prod_{i\in I} f((Y_s-Y_{(i-1)L})_{s\in \Lambda_i})$ with
\be\label{measurecost}
\E^Y_0[f((Y_s)_{0\leq s\leq L})^{-\frac{\gamma}{1-\gamma}}]\leq 2.
\ee
To make $\E^Y_0\big[\big(\Zi^{z, I}_{t,Y}\big)^\gamma\big]$ small, $f$ should be chosen to make $\E^Y_0\big[f_I(Y) \Zi^{z,I}_{t,Y}\big]$ small. There have
been two approaches in bounding $\E^Y_0\big[f_I(Y) \Zi^{z,I}_{t,Y}\big]$ in the literature.

The first approach is to choose $f_I(Y)$ to be a probability density so that $\E^Y_0\big[f_I(Y) \Zi^{z,I}_{t,Y}\big]$ becomes the annealed partition
function of a RWPM with a new law for the disorder $Y$. This approach was used in \cite{DGLT09} to prove disorder relevance for
the RPM, where the laws of the disorder at different time points are independently tilted to favor delocalization. It was later
adapted to the RWPM in dimensions $d\geq 4$ in \cite{BS09}, where the change of measure for $Y$ increases its jump rate,
which turns out to favor delocalization. To prove disorder relevance for the RPM at the critical dimension, the so-called {\em marginal disorder relevance},
which borderlines the known disorder relevance/irrelevance regimes, a more sophisticated change of measure was introduced in \cite{GLT08} for the RPM
with Gaussian disorder, which induces negative correlation between the disorder at different time points, and Gaussian calculations are used to estimate
the annealed partition function under the new disorder. For the RWPM in the critical dimension $d=3$, the analogue would be to
introduce correlation between the increments of $Y$ at different time steps. However the presence of correlation makes it unfeasible to estimate the annealed
partition function under the new disorder.

A variant approach to estimate $\E^Y_0\big[f_I(Y) \Zi^{z,I}_{t,Y}\big]$ was then introduced in \cite{L09} for the directed polymer model in random environment,
and later in \cite{GLT09} for the RPM at the critical dimension with general disorder. The function $f_I$ will be taken to be a test function on the disorder $Y$ instead
of as a probability density that changes the law of $Y$. For simplicity, $f$ in (\ref{measurecost}) is taken to be of the form
\be\label{fform}
f((Y_j)_{0\leq s\leq L}) = 1_{\{H_L(Y) \leq M\}} + \eps_M\, 1_{\{ H_L(Y) > M\}},
\ee
where $H_L(Y)$ is a functional of the disorder $Y$, positively correlated with $\Zi^{z, I}_{t,Y}$, and we choose
\be\label{epsM}
\eps_M = \P^Y_0(H_L(Y) >M)^{\frac{1-\gamma}{\gamma}}
\ee
to guarantee that (\ref{measurecost}) holds. We will make $\eps_M$ small by choosing $M$ large. To bound $\E^Y_0\big[f_I(Y) \Zi^{z,I}_{t,Y}\big]$, we use the representation (\ref{ctsZzINY}) to write
\begin{eqnarray*}
&& \!\!\!\!\!\!\!\!\!\!\!\!f_I(Y) \Zi^{z,I}_{t,Y}  \\
&=&
\!\!\!\!\!\!\!\!\!\!\!\!\!\! \int\limits_{a_1<b_1 \atop a_1,b_1\in \Lambda_{i_1}}\!\!\!\!\!\!\cdots\!\!\!\!\! \int\limits_{a_k< b_k\atop a_k,b_k\in \Lambda_{i_k}} \!\!\!\!\!\!
\prod_{j=1}^k  K(a_j-b_{j-1}) z W(a_j-b_{j-1}, Y_{a_j}-Y_{b_{j-1}}) \Zi^{z,\rm pin}_{[a_j,b_j],Y} f((Y_s-Y_{(i_j-1)L})_{s\in \Lambda_{i_j}})
\! \prod_{j=1}^k  \! {\rm d}a_j\, {\rm d}b_j,
\end{eqnarray*}
where $b_0:=0$. By Lemma \ref{L:LCLT}, the local central limit theorem  for $X$ and $X-Y$, there exists $C>0$ such that uniformly in $t>0$ and $Y$, we have
\be\label{CW}
W(t, Y_t-Y_0) = \frac{p^X_t(Y_t-Y_0)}{p^{X-Y}_t(0)} = \frac{p_t(Y_t-Y_0)}{p_{(1+\rho)t}(0)} \leq C.
\ee
Therefore
\begin{eqnarray}
&& \!\!\!\!\!\!\!\!\!\!\!\!\!\!\!\!\!\!\!\!\!\!\!\!\!\! \E^Y_0\big[f_I(Y) \Zi^{z,I}_{t,Y}\big]  \nn \\
&\leq&  \label{ctsfIZINY}
\!\!\!\!\!\!\!\!\!\!\!\! \int\limits_{a_1<b_1 \atop a_1,b_1\in \Lambda_{i_1}}\!\!\!\!\!\cdots\!\!\!\!\! \int\limits_{a_k< b_k\atop a_k,b_k\in \Lambda_{i_k}} \!\!\!\!\!\!\!\!
(Cz)^k \prod_{j=1}^k K(a_j-b_{j-1}) \E^Y_0\big[\Zi^{z,\rm pin}_{[a_j, b_j],Y} f((Y_s-Y_{(i_j-1)L})_{s\in \Lambda_{i_j}})\big] \prod_{j=1}^k {\rm d}a_j\, {\rm d}b_j,
\end{eqnarray}
where we used the independence of $(Y_s-Y_{(i-1)L})_{s\in \Lambda_i}$, $i\in\N$. The proof of (\ref{Zgammabd}), and hence Theorem \ref{T:cptctsst}, then hinges on the following key proposition.

\bp\label{P:cgZab} Let $\rho>0$, and let $L=e^{B_1\rho^{-\zeta}}$ for any fixed $\zeta>2$. Then for every $\eps>0$ and $\delta>0$, we can find suitable choices of $H_L(\cdot)$ and $M=M(L)$ in (\ref{fform}), such that for $B_1=B_1(\rho)$ sufficiently large, which can be chosen uniformly for $\rho\in(0,1]$, and for all $z\in (1, 1+L^{-1}]$, $a \in [0, (1-3\eps)L]$ and $c>L$, we have
\begin{eqnarray}
\int\limits_{a+\eps L}^{(1-\eps) L} \E^Y_0\big[\Zi^{z,\rm pin}_{[a, b],Y} f((Y_s)_{s\in [0,L]})\big] db &\leq& \delta \int\limits_a^L P(b-a) db,  \label{cgZab1}\\
\int\limits_{a+\eps L}^{(1-\eps) L} \E^Y_0\big[\Zi^{z,\rm pin}_{[a, b],Y} f((Y_s)_{s\in [0,L]})\big] K(c-b) db &\leq& \delta \int\limits_a^L P(b-a) K(c-b) db, \label{cgZab2}
\end{eqnarray}
where
\be\label{Pt1}
P(t) = \sum_{m=1}^\infty\ \idotsint\limits_{\sigma_0=0<\sigma_1<\cdots <\sigma_m=t} \prod_{i=1}^m K(\sigma_i-\sigma_{i-1}) \prod_{i=1}^{m-1} d\sigma_i,
\ee
with term for $m=1$ defined to be $K(t)$, is the renewal density associated with $K(\cdot)$.
\ep
We will show that Prop.~\ref{P:cgZab} implies the following:

\bp\label{P:coarsegrain} Let $\rho, B_1, \zeta, L, H_L$ and $M(L)$ be as in Prop.~\ref{P:cgZab}. Then for every $\eta>0$, we can choose $B_1(\rho)$ sufficiently
large, which can be chosen uniformly for $\rho\in (0,1]$, such that for all $z\in (1, 1+L^{-1}]$, $m\in\N$, and $I=\{i_1<i_2<\cdots <i_k\}\subset \{1,\cdots, m\}$, we have
\be\label{fIZINY2}
\E^Y_0\big[f_I(Y) \Zi^{z,I}_{t,Y}\big] \leq C_L \prod_{j=1}^{k} \frac{\eta}{(i_j-i_{j-1})^{\frac{3}{2}}}
\ee
for some $C_L>1$ depending only on $L$.
\ep
By (\ref{Zgsplit}) and (\ref{Holder}), Prop.~\ref{P:coarsegrain} implies that uniformly in $t=mL$, $m\in\N$, we have
\begin{eqnarray*}
\E^Y_0\big[(\Zi^z_{t,Y})^\gamma\big] \leq \sum_{k=0}^\infty \sum_{I\subset \N \atop |I|=k} \E^Y_0\big[(\Zi^{z,I}_{t,Y})^\gamma\big]
\leq \sum_{k=0}^\infty \sum_{I=\{i_1<\cdots<i_k\}} \!\!\!\!\!\!  C_L^{\gamma} \prod_{j=1}^{k} \frac{\eta^\gamma 2^{(1-\gamma)}}{(i_j-i_{j-1})^{\frac{3\gamma}{2}}}
\leq C_L^\gamma \sum_{k=0}^\infty \Big(\sum_{n=1}^\infty \frac{\eta^\gamma 2^{1-\gamma}}{n^{\frac{3\gamma}{2}}}\Big)^k,
\end{eqnarray*}
which is finite if we choose $\gamma\in (2/3,1)$, and $\eta>0$ sufficiently small such that $\sum_{n=1}^\infty \frac{\eta^\gamma 2^{1-\gamma}}{n^{\frac{3\gamma}{2}}}<1$. By the monotonicity of $\Zi^z_{t,Y}$ in $t$, this implies (\ref{Zgammabd}), and hence Theorem \ref{T:cptctsst} by the discussion following (\ref{Zgammabd}).

The key is therefore Prop.~\ref{P:cgZab}, which is the model dependent part and whose proof will be the focus of the rest of this paper.
The new idea developed  in \cite{L09} and \cite{GLT09} to bound quantities like $\E^Y_0[\Zi^{z,\rm pin}_{[a, b],Y} f((Y_s)_{s\in [0,L]})]$ is to use the renewal
representation (\ref{ctsz}) to write
\begin{eqnarray}
&&   \E^Y_0\Big[\Zi^{z,\rm pin}_{[a, b],Y} f((Y_s)_{s\in [0,L]})\Big] \nn \\
&=& \!\!
\sum_{k=1}^\infty \ \idotsint\limits_{\sigma_0=a<\cdots <\sigma_k=b} \!\!\!\!\! z^k \prod_{i=1}^k K(\sigma_i-\sigma_{i-1})
\E^Y_0\Big[f((Y_s)_{s\in [0,L]}) \prod_{i=1}^k W(\sigma_{i}-\sigma_{i-1}, Y_{\sigma_{i}}-Y_{\sigma_{i-1}})\Big]\! \prod_{i=1}^{k-1}{\rm d}\sigma_i
\nn \\
&=& \!\!
\sum_{k=1}^\infty \ \idotsint\limits_{\sigma_0=a<\cdots <\sigma_k=b} \!\!\!\!\! z^k  \prod_{i=1}^k K(\sigma_i-\sigma_{i-1})
\E^{Y^\sigma}_0[ f((Y^\sigma_s)_{s\in [0,L]})]  \prod_{i=1}^{k-1}{\rm d}\sigma_i, \label{ctsfIZINY3}
\end{eqnarray}
where $\prod_{i=1}^k W(\sigma_{i}-\sigma_{i-1}, Y_{\sigma_{i}}-Y_{\sigma_{i-1}})$ has been interpreted as the density for a change of measure for $Y$, and
$\E^{Y^\sigma}_0[\cdot]$ denotes expectation with respect to a random path $(Y^{\sigma}_s)_{0\leq s\leq t}$ whose law is absolutely continuous
with respect to that of $(Y_s)_{0\leq s\leq t}$ with density $\prod_{i=1}^k W(\sigma_i-\sigma_{i-1}, Y_{\sigma_{i}}-Y_{\sigma_{i-1}})$. Recall the form of $f$ in (\ref{fform}):
the key point is to choose the functional $H_L$ such that for typical realizations of $\sigma$ and $Y^\sigma$, $H_L((Y^\sigma_s)_{s\in [0,L]})$
is much larger than typical values of $H_L\big((Y_s)_{s\in [0,L]}\big)$. Then in (\ref{fform}), we can choose $M$ large such that $\eps_M<\!< 1$ and
$\E^{Y^\sigma}_0[ f((Y^\sigma_s)_{s\in [0,1]})]<\!< 1$. The factor $z^k$ in (\ref{ctsfIZINY3}) can be bounded by a constant of order $1$
if $z\in (1, 1+L^{-1}]$, since conditioned on the renewal process $\sigma$ with $a<b\in\sigma$, the number of renewal returns in $[a,b]$ is typically
of the order $\sqrt{b-a} \leq \sqrt{L}$.
\bigskip

\noindent
{\bf Remark.}
The above procedure applies to general weighted renewal processes in random environments, whose partition functions can be represented in the form of (\ref{ctsz})
and (\ref{ctszfree}), where given a random environment $(\Omega_s)_{s\geq 0}$ with stationary independent increments and a renewal configuration
$\sigma:=\{\sigma_0=0<\sigma_1<\cdots\}$, each two consecutive renewal times $\sigma_i<\sigma_{i+1}$ give rise to a weight factor $z W(\sigma_{i+1}-\sigma_i, (\Omega_s-\Omega_{\sigma_i})_{\sigma_i< s \leq \sigma_{i+1}})$. See e.g.\ \cite[Section 1.3]{BS09} for a more detailed exposition on how random pinning, random
walk pinning, and copolymer models can all be seen as renewal processes in random environments with different weight factors $W$. With proper normalization, $W(\sigma_{i+1}-\sigma_i, (\Omega_s-\Omega_{\sigma_i})_{\sigma_i< s \leq \sigma_{i+1}})$ can always be interpreted as a change of measure for the disorder $\Omega$.

\section{Mean and variance of $H_L(Y)$ and $H_L(Y^\sigma)$}\label{S:HLY}
We will now choose the functional $H_L(\cdot)$ in (\ref{fform}), state its essential properties, and briefly outline how these properties may lead to
Prop.~\ref{P:cgZab}. Given a renewal configuration $\sigma:=\{\sigma_0=a<\cdots<\sigma_k=b\}$, the new disorder random walk $Y^\sigma$ introduced in (\ref{ctsfIZINY3})
has heuristically smaller fluctuations than $Y$ due to the density $\prod W(\sigma_i-\sigma_{i-1}, Y_{\sigma_i}-Y_{\sigma_{i-1}}) = \prod \frac{p^X(Y_{\sigma_i}-Y_{\sigma_{i-1}})}{p^{X-Y}(\sigma_i-\sigma_{i-1})}$ which favors values of $Y_{\sigma_i}$ that are close to $Y_{\sigma_{i-1}}$.
A natural choice for $H_L$ in (\ref{fform}) is then the following. Fix $A_1:=e < A_2 <\infty$, where later in the proof of Prop.~\ref{P:cgZab} we will set $A_2=L^{\frac{1}{8}}$. The reason for this choice of $A_2$ will be explained briefly at the end of this section. Given $\zeta>2$ as in Theorem \ref{T:cptcts}, let
\be\label{xizeta}
\xi:=1-\zeta^{-1}>\frac{1}{2}.
\ee
Then we define
\be\label{HLY}
H_L(Y) = H_L((Y_s-Y_0)_{0\leq s\leq L})  : = \iint\limits_{0<r<s<L \atop A_1<s-r<A_2} \frac{1_{\{Y_r=Y_s\}}}{(\log (s-r))^\xi} drds.
\ee
We have the following bound on the mean and variance of $H_L(Y)$.

\bl\label{L:HLYbound} Let $H_L(Y)$ be defined as in (\ref{HLY}). Then
\be\label{ctsEHN}
\E^Y_0[H_L(Y)] =\iint\limits_{0<r<s<L \atop A_1<s-r<A_2} \frac{p_{\rho(s-r)}(0)}{(\log (s-r))^\xi} drds \leq (A_2-A_1)L,
\ee
and there exists some $0<C<\infty$ such that uniformly for all $A_1=e<A_2<\infty$ and $\rho>0$,
\be\label{ctsvarHN}
{\rm Var}(H_L(Y)) \leq C \rho^{-3} L.
\ee
\el
{\bf Remark.} The condition $\xi>\frac{1}{2}$ guarantees the validity of (\ref{ctsvarHN}).
The technical reason for the relation $\xi = 1-1/\zeta$, hence $\zeta>2$, will become evident in the proof Lemma~\ref{L:Hintest}, see especially (\ref{Hintest5}), below.
\bigskip

To show that in (\ref{ctsfIZINY3}), $\E^{Y^\sigma}_0[ f((Y^\sigma_s)_{s\in [0,L]})]$ is small for typical realizations of $\sigma$, it then suffices to
show that for typical realizations of $\sigma$ and $(Y^\sigma_s)_{s\in [0,L]}$, $H_L(Y^\sigma)>\E^Y_0[H_L(Y)] + D \rho^{-\frac{3}{2}}\sqrt{L}$,
where $D$ can be made arbitrarily large by choosing $B_1$ large in $L=e^{B_1\rho^{-\zeta}}$, with $B_1$ uniform for $\rho\in (0,1]$. Thus we need to
bound the mean and variance of $H_L(Y^\sigma)$ conditioned on $\sigma$. Recall that given $\sigma=\{\sigma_0=a<\sigma_1<\cdots<\sigma_k=b\}\subset [0,L]$,
the law of $(Y^\sigma_s)_{s\geq 0}$ is absolutely continuous with respect to the
law of $(Y_s)_{s\geq 0}$ with density $\prod_{i=1}^{k}
\frac{p^X_{\sigma_{i}-\sigma_{i-1}}(Y_{\sigma_{i}}-Y_{\sigma_{i-1}})}{p^{X-Y}_{\sigma_{i}-\sigma_{i-1}}(0)}$.
Due to the dependency structure of $Y^\sigma$, we will decompose $H_L(Y^\sigma)$ in (\ref{HLY}) according to whether or not the variables of integration $r<s$ satisfy $(r,s)\cap\{\sigma_0=a<\cdots<\sigma_k=b\}=\emptyset$ in order to extract some independence. Namely,
\begin{eqnarray}
\! H_L(Y^\sigma)\! =\! H_{[0,a]}^{\rm int}(Y^\sigma) + \! H_{[b,L]}^{\rm int}(Y^\sigma) + \!\! \sum_{i=1}^{k}\! H_{[\sigma_{i-1},\sigma_{i}]}^{\rm int}(Y^\sigma) + \! H_{[0,a]}^{\rm ext}(Y^\sigma) + \!\! \sum_{i=1}^{k}\! H_{[\sigma_{i-1},\sigma_{i}]}^{\rm ext}(Y^\sigma) - C_{\sigma,Y^\sigma}, \label{HLYdecomp}
\end{eqnarray}
where for any $s<t$,
\be\label{Hsigma}
\begin{aligned}
H_{[s, t]}^{\rm int}(Y^\sigma) & := \iint\limits_{s<s_1<s_2<t \atop A_1<s_2-s_1<A_2} \frac{1_{\{ Y^\sigma_{s_1}=Y^\sigma_{s_2} \}}}{(\log(s_2-s_1))^\xi} ds_2 ds_1, \\
H_{[s, t]}^{\rm ext}(Y^\sigma) & := \iint\limits_{s<s_1<t<s_2 \atop A_1<s_2-s_1<A_2} \frac{1_{\{ Y^\sigma_{s_1}=Y^\sigma_{s_2} \}}}{(\log(s_2-s_1))^\xi} ds_2 ds_1,
\end{aligned}
\ee
and
$$
C_{\sigma,Y^\sigma} := \iint\limits_{0<s_1<\sigma_k, s_2>L \atop A_1<s_2-s_1<A_2} \frac{1_{\{ Y^\sigma_{s_1}=Y^\sigma_{s_2} \}}}{(\log(s_2-s_1))^\xi} ds_2 ds_1
$$
arises because $H^{\rm ext}_{[\sigma_{i-1}, \sigma_i]}$ may include pair correlation terms $1_{\{ Y^\sigma_{s_1}=Y^\sigma_{s_2} \}}$ with $s_1<L<s_2$,
which is excluded in the definition of $H_L$. Note that
\be\label{Hextbd}
C_{\sigma,Y^\sigma} \leq A_2^2 \qquad \mbox{and} \qquad H_{[s, t]}^{\rm ext}(Y^\sigma) \leq A_2^2 \quad \quad \mbox{for all } \sigma, Y^\sigma, s<t \ \mbox{and } \rho>0 .
\ee

Conditional on $\sigma$, for any two consecutive renewal times $\sigma_{i-1}<\sigma_i$, we then have the following bounds on the mean of
$H^{\rm ext}_{[\sigma_{i-1},\sigma_i]}(Y^\sigma)$ and $H_{[\sigma_{i-1},\sigma_i]}^{\rm int}(Y^\sigma)$, and the variance of
$H_{[\sigma_{i-1},\sigma_i]}^{\rm int}(Y^\sigma)$.

\bl\label{L:Hsigmabd} For any $A_1:=e <A_2<\infty$ and $\sigma:=\{\sigma_0=a<\sigma_1<\cdots \sigma_k=b\}\subset (0,L]$, we have
\begin{gather}
\E^{Y^\sigma}_0[H_{[\sigma_{i-1}, \sigma_{i}]}^{\rm ext}(Y^\sigma)] - \E^{Y}_0[H_{[\sigma_{i-1}, \sigma_{i}]}^{\rm ext}(Y)] > 0, \qquad 1\leq i\leq k, \label{Hsigmext} \\
\E^{Y^\sigma}_0[H_{[\sigma_{i-1}, \sigma_{i}]}^{\rm int}(Y^\sigma)] - \E^{Y}_0[H_{[\sigma_{i-1}, \sigma_{i}]}^{\rm int}(Y)] > 0, \qquad 1\leq i\leq k. \label{Hsigmext2}
\end{gather}
Furthermore,
\be
\E^{Y^\sigma}_0[H_{[\sigma_0, \sigma_1]}^{\rm int}(Y^\sigma)] - \E^Y_0[H_{[\sigma_0, \sigma_1]}^{\rm int}(Y)] > \frac{C \sqrt{\sigma_1-\sigma_0}}{\sqrt{\rho} (\log(\sigma_1-\sigma_0))^\xi}\, 1_{\{2A_1<\sigma_1-\sigma_0<A_2\}} \label{Hsigmint}
\ee
and
\be
{\rm Var}(H_{[\sigma_0, \sigma_1]}^{\rm int}(Y^\sigma)|\sigma) \leq C\rho^{-3}(\sigma_1-\sigma_0), \label{Hsigmavar}
\ee
where ${\rm Var}(\cdot|\sigma)$ denotes variance w.r.t.\ $Y^\sigma$ conditioned on $\sigma$, and the $C$s in (\ref{Hsigmint})--(\ref{Hsigmavar}) are
uniform in $A_2$ and $\rho>0$.
\el

Let us sketch briefly how Lemma \ref{L:Hsigmabd} can be used to deduce Prop.~\ref{P:cgZab}. Let $\sigma$ be a renewal process conditioned on $\sigma_0=a$, and
let $Y^\sigma$ be defined accordingly by changing the measure of $Y$ independently on each renewal interval $(\sigma_{i-1}, \sigma_i)$. By the discussion following (\ref{ctsfIZINY3}), the key to proving Prop.~\ref{P:cgZab} is to show that for typical $\sigma$ and $Y^\sigma$, $H_L(Y^\sigma)$ is much larger than typical values of $H_L(Y)$. Using the decomposition (\ref{HLYdecomp}), this is achieved by controlling the mean and variance of $\sum_{i=1}^{k}\! H_{[\sigma_{i-1},\sigma_{i}]}^{\rm int}(Y^\sigma)$ and $\sum_{i=1}^{k}\! H_{[\sigma_{i-1},\sigma_{i}]}^{\rm ext}(Y^\sigma)$ conditional on $\sigma$, which is facilitated by Lemma \ref{L:Hsigmabd}.
If we only want to establish disorder relevance for a fixed $\rho>0$, which amounts to proving Prop.~\ref{P:cgZab} (more precisely, its reduction, Prop.~\ref{P:cgZfree})
by choosing $L$ sufficiently large for a fixed $\rho$, we can avoid quantitative estimates by simply applying the law of large numbers to the i.i.d.\ sequence $\big(\E^{Y^\sigma}_0[H_{[\sigma_{i-1}, \sigma_{i}]}^{\rm int}(Y^\sigma)]\big)_{i\in\N}$ and the ergodic theorem to the ergodic sequence $\big(H_{[\sigma_{i-1}, \sigma_{i}]}^{\rm ext}(Y^\sigma)\big)_{i\in\N}$ to show that,
if $L$ is sufficiently large, then for typical $\sigma$, the conditional mean of
\be\label{condmean}
\sum_{i=1}^k\! H_{[\sigma_{i-1},\sigma_{i}]}^{\rm int}(Y^\sigma)-\sum_{i=1}^k\! \E^Y_0[H_{[\sigma_{i-1},\sigma_{i}]}^{\rm int}(Y)]
\ee
far exceeds its conditional standard deviation as well as the standard deviation of $H_L(Y)$, which are of the order $C\rho^{-\frac{3}{2}}\sqrt{L}$ by
(\ref{ctsvarHN}) and (\ref{Hsigmavar}); while for typical $\sigma$ and $Y^\sigma$,
$$
\sum_{i=1}^k\! H_{[\sigma_{i-1},\sigma_{i}]}^{\rm ext}(Y^\sigma)- \sum_{i=1}^k\! \E^Y_0[H_{[\sigma_{i-1},\sigma_{i}]}^{\rm ext}(Y)] >0.
$$
To get quantitative bounds on the gap between the annealed and the quenched critical points, we need to get bounds on $L$, and this is achieved by replacing the
law of large numbers above with a quantitative estimate on $\sum_{i=1}^k\! \E^{Y^\sigma}_0[H_{[\sigma_{i-1},\sigma_{i}]}^{\rm int}(Y^\sigma)]$, and replacing the ergodic
theorem with a quantitative bound on the conditional variance of $\sum_{i=1}^k\! H_{[\sigma_{i-1},\sigma_{i}]}^{\rm ext}(Y^\sigma)$ and then applying the Markov inequality. The details will be given in Sec.~\ref{S:cgZab}.

The reason for choosing $A_2=L^{\frac{1}{8}}$ in the definition of $H_L(Y)$ is the following. When
we lower bound the conditional mean in (\ref{condmean}) (conditional on $\sigma$) using (\ref{Hsigmint}), we need to choose $A_2$ as large as possible. It turns out that any power of $L$
will suffice, as will be seen in the proof of Lemma~\ref{L:sumZiL}. On the other hand, when we upper bound the conditional variance of $\sum_{i=1}^k\! H_{[\sigma_{i-1},\sigma_{i}]}^{\rm ext}(Y^\sigma)$, we need to choose $A_2$ to be a sufficiently small power of $L$, which can be seen from the bounds in
(\ref{Hextbd}) as well as in the proof of Lemma \ref{L:Hextest}. The choice $A_2=L^{\frac{1}{8}}$ turns out to be sufficient for our purposes.

\section{Proof of Proposition \ref{P:cgZab}}\label{S:cgZab}

We now prove Prop.~\ref{P:cgZab} using the functional $H_L$ defined in (\ref{HLY}) and Lemmas \ref{L:HLYbound} and \ref{L:Hsigmabd}. We remark that Prop.~\ref{P:cgZab}
is the analogue of \cite[Lemma 3.1]{GLT09} formulated for the discrete time random pinning model. The main difference is that \cite[Lemma 3.1]{GLT09} involves a comparison of the integrands on both sides of (\ref{cgZab2}) for each $b$ with $b-a\geq \eps L$. Our integral formulation of (\ref{cgZab2}) allows us to reduce more easily estimates involving renewal configurations pinned at two points $a<b$ to renewal configurations pinned only at $a$. More precisely, we reduce Prop.~\ref{P:cgZab}
to the following proposition by extracting some model independent renewal calculations.

\bp\label{P:cgZfree}
Let $A_1=e$ and $A_2=L^{\frac{1}{8}}$ in the definition of $H_L$ in (\ref{HLY}), with $L=e^{B_1\rho^{-\zeta}}$ as in Prop.~\ref{P:cgZab}. Then for every
$\eps>0$ and $\delta>0$, we can find $D>0$ and set $M=\E^Y_0[H_L(Y)]+D\rho^{-\frac{3}{2}}\sqrt{L}$ in (\ref{fform}), such that for all $B_1=B_1(\rho)$ sufficiently
large, which can be chosen uniformly for $\rho\in (0,1]$, and for all $z\in (1, 1+L^{-1}]$ and $a\in [0, (1-3\eps)L]$, we have
\be\label{cgZfreeP}
\int\limits_a^L \E^Y_0\big[\Zi^{z,\rm pin}_{[a, b],Y} f((Y_s)_{s\in [0,L]})\big] K([L-b,\infty)) db \leq \delta,
\ee
where $K([x,\infty))=\int_x^\infty K(t)dt = \int_x^\infty (1+\rho)G^{-1}p_{(1+\rho)t}(0)dt$.
\ep
{\bf Proof of Prop.~\ref{P:cgZab}.} The deduction of Prop.~\ref{P:cgZab} from Prop.~\ref{P:cgZfree} is model independent and depends only on $K(\cdot)$. Since $K(t)\sim \frac{C}{(1+\rho)^\frac{1}{2}t^{\frac{3}{2}}}$, we have $K([t,\infty))=\int_t^\infty K(s)ds \sim \frac{2C}{\sqrt{(1+\rho)t}}$. By Lemma
\ref{L:renewalprob}, we also have $P(t)\sim \frac{C}{\sqrt{(1+\rho)t}}$, where $P(t)$ is defined in (\ref{Pt1}). Therefore, given $\eps>0$ and for $B_1$ (and hence $L$) large, there exist $C_1$ and $C_2$ depending only on $\eps>0$,
such that uniformly for all $\eps L \leq a+\eps L\leq b_1,b_2\leq (1-\eps)L$ and $c>L$, we have
\be\label{cgZabpf1}
 C_1 \leq \frac{P(b_1-a)}{P(b_2-a)} \leq C_2, \qquad C_1 \leq \frac{K(c-b_1)}{K(c-b_2)}\leq C_2.
\ee
Under the assumptions of Prop.~\ref{P:cgZab}, by Prop.~\ref{P:cgZfree}, we have
\be\nn
\delta \geq \!\!\!\! \int\limits_{a+\eps L}^{(1-\eps)L} \frac{\E^Y_0\big[\Zi^{z,\rm pin}_{[a, b],Y} f((Y_s)_{s\in [0,L]})\big]}{P(b-a)} P(b-a) K([L-b,\infty)) db
\geq \frac{C}{L} \!\!\!\! \int\limits_{a+\eps L}^{(1-\eps) L} \frac{\E^Y_0\big[\Zi^{z,\rm pin}_{[a, b],Y} f((Y_s)_{s\in [0,L]})\big]}{P(b-a)} db,
\ee
where $C$ depends only on $\eps$ and $\rho$ and is uniform for $\rho\in (0,1]$. Together with (\ref{cgZabpf1}), this implies
\begin{eqnarray*}
&& \int\limits_{a+\eps L}^{(1-\eps) L} \E^Y_0\big[\Zi^{z,\rm pin}_{[a, b],Y} f((Y_s)_{s\in [0,L]})\big] K(c-b) db \\
&=& \int\limits_{a+\eps L}^{(1-\eps) L} \frac{\E^Y_0\big[\Zi^{z,\rm pin}_{[a, b],Y} f((Y_s)_{s\in [0,L]})\big]}{P(b-a)} P(b-a) K(c-b) db \\
&=& \int\limits_{a+\eps L}^{(1-\eps) L} \frac{\E^Y_0\big[\Zi^{z,\rm pin}_{[a, b],Y} f((Y_s)_{s\in [0,L]})\big]}{P(b-a)}
\bigg(\int\limits_{a+\eps L}^{(1-\eps)L} \hspace{-0.3em}  \frac{P(b_2-a) K(c-b_2)}{P(b-a) K(c-b)} db_2 \bigg)^{-1} \hspace{-0.3em} db
\int\limits_{a+\eps L}^{(1-\eps)L} \hspace{-0.9em} P(b-a) K(c-b) db \\
&\leq& \frac{1}{C_1^2} \frac{1}{(1-\eps)L-(a+\eps L)} \frac{\delta L}{C} \int\limits_a^{L} P(b-a) K(c-b) db
\leq \frac{\delta}{\eps C C_1^2} \int\limits_a^{L} P(b-a) K(c-b) db,
\end{eqnarray*}
where we used the assumption that $a\leq (1-3\eps)L$. Since given $\eps>0$, we can choose $\delta>0$ arbitrarily small by Prop.~\ref{P:cgZfree}, (\ref{cgZab2}) then
follows. The proof of (\ref{cgZab1}) is similar. Note that we need (\ref{cgZab1}) because we are studying $\Zi^{z}_{t, Y}$, instead of the constrained partition function $\Zi^{z,\rm pin}_{t,Y}$ as in \cite{GLT09}.
\qed
\bigskip

The proof of Prop.~\ref{P:cgZfree} is based on two lemmas, for which we first introduce some notation. Let $\sigma:=\{\sigma_0=a<\sigma_1<\cdots\}$ be a renewal process on $[a,\infty)$ with renewal time distribution $K(t)dt=(1+\rho)G^{-1}p_{(1+\rho)t}(0)dt$, and let $k(\sigma, L):=|\sigma\cap (a,L]|$. Let $\P^\sigma_a$ and $\E^\sigma_a$
denote respectively probability and expectation for $\sigma$. Let $(\tilde Y^\sigma_s)_{s\geq 0}$ be defined analogously to $Y^\sigma$, where conditional on $\sigma$, for each $n\in\N$, the law of $(\tilde Y^\sigma_s)_{0\leq s\leq \sigma_n}$ is absolutely continuous w.r.t.\ the law of $(Y_s)_{0\leq s\leq \sigma_n}$ with Radon-Nikodym derivative
$$
\prod_{i=1}^n \frac{p_{\sigma_i-\sigma_{i-1}}(Y_{\sigma_i}-Y_{\sigma_{i-1}})}{p_{(1+\rho)(\sigma_i-\sigma_{i-1})}(0)}.
$$
Then we have

\bl\label{L:Hintest} Let $L=e^{B_1\rho^{-\zeta}}$, $A_1=e$ and $A_2=L^{\frac{1}{8}}$.  Let $a\in [0, (1-\eps)L]$ for some $0<\eps<1$. For any
$D_1>0$ and $\delta>0$, if $B_1=B_1(\rho)$ is sufficiently large, which can be chosen uniformly for $\rho\in (0,1]$, then for all $a\in [0, (1-\eps)L]$, we have
\be\label{Hintest1}
\P^\sigma_a\P^{\tilde Y^\sigma}_0\Big(\sum_{i=1}^{k(\sigma, L)}H^{\rm int}_{[\sigma_{i-1}, \sigma_i]}(\tilde Y^\sigma) - \sum_{i=1}^{k(\sigma,L)}\E^Y_0[H^{\rm int}_{[\sigma_{i-1},\sigma_i]}(Y)] < D_1\rho^{-\frac{3}{2}}\sqrt{\eps L} \Big) \leq \delta.
\ee
\el

\bl\label{L:Hextest} Let $L$, $A_1$, $A_2$, $\eps$ and $a$ be as in Lemma \ref{L:Hintest}. For any $\delta>0$, if $B_1=B_1(\rho)$ is sufficiently large, which
can be chosen uniformly for $\rho\in (0,1]$, then for all $a\in [0, (1-\eps)L]$, we have
\be\label{Hextest1}
\P^\sigma_a\P^{\tilde Y^\sigma}_0\Big(\sum_{i=1}^{k(\sigma, L)}H^{\rm ext}_{[\sigma_{i-1}, \sigma_i]}(\tilde Y^\sigma) - \sum_{i=1}^{k(\sigma,L)}\E^Y_0[H^{\rm ext}_{[\sigma_{i-1},\sigma_i]}(Y)] < -\sqrt{L} \Big) \leq \delta.
\ee
\el
We defer the proofs of Lemmas \ref{L:Hintest}--\ref{L:Hextest} and first deduce Prop.~\ref{P:cgZfree}.
\bigskip

\noindent
{\bf Proof of Prop.~\ref{P:cgZfree}.} Let $\sigma$, $k:=k(\sigma, L)$ and $\tilde Y^\sigma$ be as introduced before Lemma \ref{L:Hintest}. Then we define
$Y^\sigma_s=\tilde Y^\sigma_s$ for $s\in [0, \sigma_{k}]$, and $Y^\sigma_s-Y^\sigma_{\sigma_{k}}=Y_s-Y_{\sigma_{k}}$ for $s\geq \sigma_{k}$.
By (\ref{ctsfIZINY3}), uniformly in $z\in (1, 1+1/L]$, we may rewrite the LHS of (\ref{cgZfreeP}) as
\begin{eqnarray}
&& \int\limits_a^L \E^Y_0\big[\Zi^{z,\rm pin}_{[a, b],Y} f((Y_s)_{s\in [0,L]})\big] K([L-b,\infty)) db
= \E^\sigma_a\Big[1_{\{k\geq 1\}}z^{k}\E^{Y^\sigma}_0\big[f((Y^\sigma_s)_{s\in [0,L]})\big]\Big] \nn \\
&\leq& \E^\sigma_a[(1+L^{-1})^{k} 1_{\{k>L\}}] + e\E^\sigma_a\big[\E^{Y^\sigma}_0\big[f((Y^\sigma_s)_{s\in [0,L]})\big]\big] \nn \\
&\leq& C_1 e^{-C_2L} + \ e\P^\sigma_a\P^{Y^{\sigma}}_0(H_L(Y^\sigma)\leq M) + e\eps_M, \label{cgZfree0}
\end{eqnarray}
where the bound for $\E^\sigma_a[(1+L^{-1})^{k} 1_{\{k>L\}}]$ follows from standard large deviation estimates for the i.i.d.\ random
variables $(\sigma_i-\sigma_{i-1})_{i\in\N}$. In (\ref{cgZfree0}), the first and the last terms can both be made arbitrarily small
by choosing $B_1$ large enough, and $D$ large enough in $M=\E^Y_0[H_L(Y)]+D\rho^{-\frac{3}{2}}\sqrt{L}$, which
follows from (\ref{epsM}) and (\ref{ctsvarHN}).

Recall the decomposition of $H_L(Y^\sigma)$ in (\ref{HLYdecomp}). Fix any $\delta>0$. Since $\tilde Y^\sigma_s=Y^\sigma_s$ on $[0, \sigma_k]$, by Lemma \ref{L:Hintest},
we can choose $B_1$ large enough (uniformly for $\rho\in (0,1]$) such that for all $a\in [0, (1-3\eps)L]$, we have
\be\label{cgZfree1}
\P^\sigma_a\P^{Y^\sigma}_0\Big(\sum_{i=1}^{k(\sigma, L)}H^{\rm int}_{[\sigma_{i-1}, \sigma_i]}(Y^\sigma) - \sum_{i=1}^{k(\sigma,L)}\E^Y_0[H^{\rm int}_{[\sigma_{i-1},\sigma_i]}(Y)] < 4D\rho^{-\frac{3}{2}}\sqrt{L} \Big) \leq \delta.
\ee
By the same reasoning as in (\ref{Hextbd}), we note that $\sum_{i=1}^{k(\sigma,L)} H_{[\sigma_{i-1}, \sigma_{i}]}^{\rm ext}(\tilde Y^\sigma)$ and $\sum_{i=1}^{k(\sigma,L)} H_{[\sigma_{i-1}, \sigma_{i}]}^{\rm ext}(Y^\sigma)$ differ by at most $A_2^2=L^{\frac{1}{4}}$. Therefore by Lemma \ref{L:Hextest}, we can choose $B_1$ large
enough such that
\be\label{cgZfree2}
\P^\sigma_a\P^{Y^\sigma}_0\Big(\sum_{i=1}^{k(\sigma, L)}H^{\rm ext}_{[\sigma_{i-1}, \sigma_i]}(Y^\sigma) - \sum_{i=1}^{k(\sigma,L)}\E^Y_0[H^{\rm ext}_{[\sigma_{i-1},\sigma_i]}(Y)] < -2\sqrt{L} \Big) \leq \delta.
\ee
In the decomposition of $H_L(Y^\sigma)$ in (\ref{HLYdecomp}), again by (\ref{Hextbd}), we have
\be\label{cgZfree3}
H_{[0, a]}^{\rm ext}(Y^\sigma) - \E^Y_0[H_{[0, a]}^{\rm ext}(Y)] - C_{\sigma, Y^\sigma} + \E^Y_0[C_{\sigma,Y}] \geq -2 A_2^2 = -2L^{\frac{1}{4}},
\ee
where $C_{\sigma,Y}$ is defined exactly as $C_{\sigma, Y^\sigma}$ with $Y^\sigma$ replaced by $Y$. The same calculation as in the proof of Lemma \ref{L:HLYbound} (\ref{ctsvarHN}), to appear in Section~\ref{S:HLYbound}, shows that ${\rm Var}(H_{[0,a]}^{\rm int}(Y))\leq C\rho^{-3}L$ and ${\rm Var}(H_{[\sigma_k,L]}^{\rm int}(Y))\leq C\rho^{-3}L$. Since by construction, $Y^\sigma$ have the same increments as $Y$ on $[0,a]$ and $[\sigma_k, L]$, we can choose $D$ large enough such that
\be\label{cgZfree4}
\begin{aligned}
\P^\sigma_a\P^{Y^\sigma}_0\Big(H_{[0, a]}^{\rm int}(Y^\sigma) - \E^Y_0[H_{[0, a]}^{\rm int}(Y)]  &\ \leq \  -D\rho^{-\frac{3}{2}}\sqrt{L} \Big)
\leq \delta, \\
\P^\sigma_a\P^{Y^\sigma}_0\Big(H_{[\sigma_k, L]}^{\rm int}(Y^\sigma) - \E^Y_0[H_{[\sigma_k, L]}^{\rm int}(Y)] &\ \leq\ -D\rho^{-\frac{3}{2}}\sqrt{L} \Big)
\leq \delta.
\end{aligned}
\ee
If we first choose $D$ large and then $B_1$ large, and let $M=\E^Y_0[H_L(Y)]+ D\rho^{-\frac{3}{2}}\sqrt{L}$, then by the decomposition of $H_L(Y^\sigma)$ in (\ref{HLYdecomp}) and (\ref{cgZfree1})--(\ref{cgZfree4}), we find that in (\ref{cgZfree0}), we have
$$
\P^\sigma_a\P^{Y^\sigma}_0( H_L(Y^\sigma) \leq M) \leq 4\delta.
$$
Since $\delta>0$ can be made arbitrarily small, Prop.~\ref{P:cgZfree} then follows.
\qed
\bigskip

We now prove Lemmas \ref{L:Hintest}--\ref{L:Hextest} by controlling the mean and variance of $\sum_{i=1}^{k(\sigma, L)}H^{\rm int}_{[\sigma_{i-1},\sigma_i]}(\tilde Y^\sigma)$ and $\sum_{i=1}^{k(\sigma, L)}H^{\rm ext}_{[\sigma_{i-1},\sigma_i]}(\tilde Y^\sigma)$, conditional on $\sigma$.
Our bound on the conditional mean of $\sum_{i=1}^{k(\sigma, L)}H^{\rm int}_{[\sigma_{i-1},\sigma_i]}(\tilde Y^\sigma)$ is based on the following lemma,
which also leads to our choice of $L=e^{B_1\rho^{-\zeta}}$.

\bl\label{L:sumZiL} Let $(\Delta_i)_{i\in\N}$ be i.i.d.\ with common distribution $K(t)dt$ on $[0,\infty)$, where we have $K(t)=(1+\rho)p_{(1+\rho)t}(0)/G$. Let $Z^L_i = \frac{\sqrt{\Delta_i}}{(\log\Delta_i)^\xi}1_{\{2e<\Delta_i<L^{\frac{1}{8}}\}}$ with $L= e^{B_1\rho^{-\zeta}}$, where $(1-\xi)\zeta=1$ and $\zeta>2$.
Then there exists $B_2>0$ such that for any $h>0$ and $\delta>0$, if $B_1=B_1(\rho)$ is sufficiently large, which can be chosen uniformly for $\rho\in (0,1]$,
then we have
\be\label{sumZiL}
\P\Big(\sum_{i=1}^{h\sqrt{L}} Z^L_i < B_2 h\sqrt{L}(\log L)^{1-\xi}\Big) \leq \delta.
\ee
\el
{\bf Proof.} If $B_1$ is sufficiently large, which can be chosen uniformly for $\rho\in (0,1]$, we have
$$
\mu_L := \E[Z^L_i] = \int_{2e}^{L^{\frac{1}{8}}} \frac{\sqrt{\Delta}}{(\log \Delta)^\xi} \frac{(1+\rho)p_{(1+\rho)\Delta}(0)}{G} d\Delta \geq
C \int_{2e}^{L^{\frac{1}{8}}} \frac{d\Delta}{\Delta (\log \Delta)^\xi} \geq 4B_2 (\log L)^{1-\xi}
$$
for some $B_2>0$ independent of $B_1$, $\xi\in (\frac{1}{2},1)$ and $\rho \in (0,1]$. We then prove (\ref{sumZiL}) by a large deviation estimate.
Namely, if we let $M(\lambda) = \log \E[e^{\lambda Z^L_1}]$, then for any $\lambda<0$, we have
\be\label{sumZiL2}
\P\Big(\sum_{i=1}^{h\sqrt{L}} Z^L_i < B_2 h\sqrt{L}(\log L)^{1-\xi}\Big) \leq \exp\big\{h\sqrt{L}(M(\lambda) - \lambda B_2 (\log L)^{1-\xi})\big\}.
\ee
Let $\lambda = - L^{\frac{1}{8}}$. Then for $B_1$ sufficiently large, uniformly in $\rho \in (0,1]$, we have
\begin{eqnarray*}
M(\lambda) &=& \log \Big(\int_0^\infty e^{-\frac{\sqrt{\Delta}}{L^{\frac{1}{8}}(\log \Delta)^\xi} 1_{\{2e<\Delta<L^{\frac{1}{8}}\}}} \frac{(1+\rho)p_{(1+\rho)\Delta}(0)}{G} d\Delta\Big) \\
&<& \log \Big(\int_0^\infty \big(1- \frac{\sqrt{\Delta}}{2L^{\frac{1}{8}}(\log \Delta)^\xi} 1_{\{2e<\Delta<L^{\frac{1}{8}}\}}\big)  \frac{(1+\rho)p_{(1+\rho)\Delta}(0)}{G} d\Delta\Big) \\
&=& \log\Big(1-\frac{\mu_L}{2L^{\frac{1}{8}}}\Big) < - \frac{\mu_L}{2L^{\frac{1}{8}}} \leq 2 \lambda B_2 (\log L)^{1-\xi}.
\end{eqnarray*}
Therefore the RHS of (\ref{sumZiL2}) is bounded by $\exp\{-\frac{h\sqrt{L}B_2 (\log L)^{1-\xi}}{L^{\frac{1}{8}}}\}$, which tends to $0$ uniformly in $\rho\in(0,1]$ as
$B_1\uparrow\infty$, thus implying (\ref{sumZiL}).
\qed
\bigskip

\noindent
{\bf Proof of Lemma \ref{L:Hintest}.} For $i\in\N$, let $\Delta_i=\sigma_i-\sigma_{i-1}$, which are i.i.d.\ with common distribution $K(t)dt=(1+\rho)p_{(1+\rho)t}(0)dt/G$.
By Lemma \ref{L:renewalcounts}, for any $\delta>0$, we can find $C_1>0$ small enough such that for all $L$ sufficiently large and uniformly in $a\in [0,(1-\eps)L]$,
we have
\be\label{Hintest2}
\P^\sigma_a(k(\sigma, L) <C_1\sqrt{\eps L}) \leq \frac{\delta}{4}.
\ee
By Lemma \ref{L:Hsigmabd} (\ref{Hsigmint}), almost surely with respect to $\sigma$,
\be\label{Hintest3}
\sum_{i=1}^{k(\sigma, L)}\Big(\E^{\tilde Y^\sigma}_0[H^{\rm int}_{[\sigma_{i-1}, \sigma_i]}(\tilde Y^\sigma)] - \E^Y_0[H^{\rm int}_{[\sigma_{i-1},\sigma_i]}(Y)]\Big)
> \sum_{i=1}^{k(\sigma, L)} \frac{C \sqrt{\Delta_i}}{\sqrt{\rho} (\log\Delta_i)^\xi} 1_{\{2e <\Delta_i<L^{\frac{1}{8}}\}},
\ee
while by Lemma \ref{L:sumZiL}, given $h=C_1\sqrt{\eps}$ and $\delta>0$, we can find $B_2>0$ such that for all $B_1$
sufficiently large, we have
\be\label{Hintest4}
\P^\sigma_a\Big(\sum_{i=1}^{C_1\sqrt{\eps L}} \frac{C \sqrt{\Delta_i}}{\sqrt{\rho} (\log\Delta_i)^\xi} 1_{\{2e <\Delta_i<L^{\frac{1}{8}}\}} < \frac{C B_2}{\sqrt{\rho}} C_1\sqrt{\eps L}(\log L)^{1-\xi}\Big) \leq \frac{\delta}{4}.
\ee
Therefore for a set of $\sigma$ with probability at least $1-\frac{\delta}{2}$, conditional on $\sigma$, we have
\be\label{Hintest5}
\sum_{i=1}^{k(\sigma, L)}\Big(\E^{\tilde Y^\sigma}_0[H^{\rm int}_{[\sigma_{i-1}, \sigma_i]}(\tilde Y^\sigma)] - \E^Y_0[H^{\rm int}_{[\sigma_{i-1},\sigma_i]}(Y)]\Big)
\geq CC_1B_2 \frac{\sqrt{\eps L}}{\sqrt\rho}(\log L)^{1-\xi} = D_2\rho^{-\frac{3}{2}}\sqrt{\eps L},
\ee
where $D_2= CC_1B_2B_1^{1-\xi}$ and we used $L=e^{B_1\rho^{-\zeta}}=e^{B_1\rho^{\frac{-1}{1-\xi}}}$. For any $\sigma$, by Lemma \ref{L:Hsigmabd} (\ref{Hsigmavar}),
$$
{\rm Var}\Big(\sum_{i=1}^{k(\sigma, L)}H^{\rm int}_{[\sigma_{i-1}, \sigma_i]}(\tilde Y^\sigma)\Big|\sigma\Big) \leq C\rho^{-3}(L-a)\leq C\rho^{-3}L.
$$
$D_2\sqrt{\eps}$ can be made arbitrarily large by choosing $B_1$ large. Therefore for any $D_1>0$ and for all $\sigma$ satisfying (\ref{Hintest5}), by making $B_1$
sufficiently large, we have by Markov inequality
$$
\P^{\tilde Y^\sigma}_0\Big( \sum_{i=1}^{k(\sigma, L)}H^{\rm int}_{[\sigma_{i-1}, \sigma_i]}(\tilde Y^\sigma) - \sum_{i=1}^{k(\sigma,L)}\E^Y_0[H^{\rm int}_{[\sigma_{i-1},\sigma_i]}(Y)] < D_1\rho^{-\frac{3}{2}}\sqrt{\eps L} \Big) \leq \frac{\delta}{2}.
$$
Since the set of $\sigma$ that violates (\ref{Hintest5}) has probability at most $\frac{\delta}{2}$, this implies (\ref{Hintest1}).
\qed
\bigskip

\noindent
{\bf Proof of Lemma \ref{L:Hextest}.} By Lemma \ref{L:Hsigmabd}~(\ref{Hsigmext}), for all $\sigma$, we have
$$
\sum_{i=1}^{k(\sigma, L)}\E^{\tilde Y^\sigma}_0\big[H^{\rm ext}_{[\sigma_{i-1},\sigma_i]}(\tilde Y^\sigma)\big] - \sum_{i=1}^{k(\sigma, L)}\E^Y_0\big[H^{\rm ext}_{[\sigma_{i-1},\sigma_i]}(Y)\big]>0.
$$
Therefore to establish (\ref{Hextest1}), by Markov inequality, it suffices to show that: For any $\delta>0$, if $B_1$ in $L=e^{B_1\rho^{-\zeta}}$
is sufficiently large (uniform for $\rho \in (0,1]$), then for all $a\in [0, (1-\eps)L]$, we have
\be\label{Hextest2}
\P^\sigma_a\Big( \mbox{Var}\Big(\sum_{i=1}^{k(\sigma, L)}H^{\rm ext}_{[\sigma_{i-1},\sigma_i]}(\tilde Y^\sigma) \Big| \sigma \Big)\geq \delta L\Big)<\delta.
\ee
We will decompose the sum $\sum_{i=1}^{k(\sigma,L)}H^{\rm ext}_{[\sigma_{i-1},\sigma_i]}(\tilde Y^\sigma)$ to extract some independence.

Given $\sigma:=\{\sigma_0=a<\sigma_1<\cdots\}$, let $\tau_0=0$, and for $j\in\N$, define inductively
$$
\tau_j := \min\{ i> \tau_{j-1}: \sigma_i-\sigma_{i-1}\geq A_2=L^{\frac{1}{8}}\}.
$$
Let $J=J(\sigma, L):= \max\{j\in\N: \sigma_{\tau_j}\leq L\}$. Then we have the decomposition
$$
\sum_{i=1}^{k(\sigma,L)}H^{\rm ext}_{[\sigma_{i-1},\sigma_i]}(\tilde Y^\sigma) = \sum_{j=1}^J \sum_{i=\tau_{j-1}+1}^{\tau_j-1}H^{\rm ext}_{[\sigma_{i-1},\sigma_i]}(\tilde Y^\sigma) + \sum_{i=\tau_J+1}^{k(\sigma,L)}H^{\rm ext}_{[\sigma_{i-1},\sigma_i]}(\tilde Y^\sigma) + \sum_{j=1}^J H^{\rm ext}_{[\sigma_{\tau_j-1},\sigma_{\tau_j}]}(\tilde Y^\sigma).
$$
Now note that conditional on $\sigma$, $\sum_{i=\tau_{j-1}+1}^{\tau_j-1}H^{\rm ext}_{[\sigma_{i-1},\sigma_i]}(\tilde Y^\sigma)$, $1\leq j\leq J$, and $\sum_{i=\tau_J+1}^{k(\sigma,L)}H^{\rm ext}_{[\sigma_{i-1},\sigma_i]}(\tilde Y^\sigma)$ are all independent. Similarly, for even (resp.\ odd) $1\leq j\leq J$,
$H^{\rm ext}_{[\sigma_{\tau_j-1},\sigma_{\tau_j}]}(\tilde Y^\sigma)$ are all independent. Therefore using independence and the fact that
$\mbox{Var}(X+Y+Z)\leq 3(\mbox{Var}(X)+\mbox{Var}(Y)+\mbox{Var}(Z))$
(with $Y:=\sum\limits_{j \, \mathrm{even}}^J H^{\rm ext}_{[\sigma_{\tau_j-1},\sigma_{\tau_j}]}(\tilde Y^\sigma)$, $Z:=\sum\limits_{j \, \mathrm{odd}}^J H^{\rm ext}_{[\sigma_{\tau_j-1},\sigma_{\tau_j}]}(\tilde Y^\sigma)$, $X:=\sum\limits_{i=1}^{k(\sigma,L)}H^{\rm ext}_{[\sigma_{i-1},\sigma_i]}(\tilde Y^\sigma)-Y-Z$ ), we have
\begin{eqnarray}
\!\!\!\!\!\!&&\!\!\!\!\!\! V(\sigma,L):=\mbox{Var}\Big(\sum_{i=1}^{k(\sigma, L)}H^{\rm ext}_{[\sigma_{i-1},\sigma_i]}(\tilde Y^\sigma) \Big| \sigma \Big)
\leq  3 \sum_{j=1}^J \mbox{Var}\Big(\sum_{i=\tau_{j-1}+1}^{\tau_j-1}H^{\rm ext}_{[\sigma_{i-1},\sigma_i]}(\tilde Y^\sigma)\Big|\sigma\Big) \label{Hextest3} \\
&& \qquad\qquad\qquad\qquad\qquad\quad  + 3\mbox{Var}\Big(\sum_{i=\tau_J+1}^{k(\sigma,L)}H^{\rm ext}_{[\sigma_{i-1},\sigma_i]}(\tilde Y^\sigma)\Big|\sigma\Big) + 3\sum_{j=1}^J \mbox{Var}\Big(H^{\rm ext}_{[\sigma_{\tau_j-1},\sigma_{\tau_j}]}(\tilde Y^\sigma)\Big|\sigma\Big). \nn
\end{eqnarray}
By (\ref{Hextbd}), $H^{\rm ext}_{[\sigma_{\tau_j-1},\sigma_{\tau_j}]}(\tilde Y^\sigma) \leq A_2^2$ for each $j\in\N$. Similarly, using the definition of $H^{\rm ext}$
in (\ref{Hsigma}),
$$
\sum_{i=\tau_{j-1}+1}^{\tau_j-1}H^{\rm ext}_{[\sigma_{i-1},\sigma_i]}(\tilde Y^\sigma) \leq (\sigma_{\tau_j-1}-\sigma_{\tau_{j-1}})A_2
\qquad  \mbox{for }\ j\in\N.
$$
Therefore we obtain from (\ref{Hextest3})
\be\label{Hextest4}
\E^\sigma_a[V(\sigma, L)] \leq 3 \E^\sigma_a\Big[\sum_{j=1}^{J+1} (\sigma_{\tau_j-1}-\sigma_{\tau_{j-1}})^2A_2^2\Big] + 3\E^\sigma_a[JA_2^4].
\ee
Note that $\sigma_{\tau_{J+1}-1} \geq \sigma_{k(\sigma,L)}$, in particular, the second term on the right-hand side of (\ref{Hextest3}) is accounted for in (\ref{Hextest4}).

Let $\Delta_i=\sigma_i-\sigma_{i-1}$, which are i.i.d. Then $k(\sigma, L)+1$ is a stopping time w.r.t.\ the sequence $(\Delta_i)_{i\in\N}$, and by Wald's
equation~\cite[Sec.~3.1]{D96},
\begin{eqnarray*}
&& \E^\sigma_a[J]\leq \E^\sigma_a\Big[ \sum_{i=1}^{k(\sigma, L)+1} 1_{\{\Delta_i\geq A_2\}}\Big] = \E^\sigma_a[1+k(\sigma,L)] \P^\sigma_a(\Delta_1\geq A_2) \\
&\leq& (1+C_1\sqrt{L}) \int_{A_2}^\infty \frac{(1+\rho)p_{(1+\rho)t}(0)dt}{G} \leq C\sqrt{\frac{L}{A_2}},
\end{eqnarray*}
where we used Lemma \ref{L:renewalcounts} and the discussion following it to deduce $\E^\sigma_a[k(\sigma,L)]=C_1\sqrt{L}$, and we bounded $p_t(0)$ by
$Ct^{-\frac{3}{2}}$. Note that $(\sigma_{\tau_{j}-1}-\sigma_{\tau_{j-1}}, \sigma_{\tau_j}-\sigma_{\tau_j-1})$, $j\in\N$, is an i.i.d.\ sequence of
$\R^2$-valued random variables, and $J+1$ is a stopping time with respect to this sequence. So again by Wald's equation,
$$
\E^\sigma_a\Big[\sum_{j=1}^{J+1} (\sigma_{\tau_j-1}-\sigma_{\tau_{j-1}})^2\Big] = \E^\sigma_a[J+1]\E^\sigma_a[(\sigma_{\tau_1-1}-\sigma_{\tau_0})^2]
\leq (1+C\sqrt{L/A_2}) \E^\sigma_a\Big[\big(\sum_{i=1}^{\tau_1} \Delta_i 1_{\{\Delta_i<A_2\}}\big)^2\Big],
$$
where if we denote $Z_i:=\Delta_i 1_{\{\Delta_i<A_2\}}$ and $\mu:=\E^\sigma_a[Z_i]=\int_0^{A_2}(1+\rho)G^{-1}tp_{(1+\rho)t}(0)dt\leq C\sqrt{A_2}$, then
\begin{eqnarray*}
&& \E^\sigma_a\Big[\big(\sum_{i=1}^{\tau_1} \Delta_i 1_{\{\Delta_i<A_2\}}\big)^2\Big] = \E^\sigma_a\Big[\big(\sum_{i=1}^{\tau_1} Z_i\big)^2\Big] \\
&\leq& 2\E^\sigma_a\Big[\big( \sum_{i=1}^{\tau_1}(Z_i-\mu)\big)^2\Big] + 2\mu^2\E^\sigma_a[\tau_1^2]
= 2\E^\sigma_a[\tau_1] \E^\sigma_a[(Z_1-\mu)^2] + 2\mu^2 \E^\sigma_a[\tau_1^2] \\
&\leq& 2 \E^\sigma_a[\tau_1] \E^\sigma_a[Z_1^2] + 2\mu^2 \E^\sigma_a[\tau_1^2]
\leq \frac{2\int_0^{A_2}(1+\rho)G^{-1}t^2p_{(1+\rho)t}(0)dt}{\P^\sigma_a(\Delta_1\geq A_2)} + \frac{4C^2A_2}{\P^\sigma_a(\Delta_1\geq A_2)^2} \\
&\leq& C A_2^2,
\end{eqnarray*}
where we have used Wald's second equation~\cite[Sec.~3.1]{D96} and the fact that $\tau_1$ is a stopping time for $(\Delta_i)_{i\in\N}$, which is geometrically distributed
with $\E^\sigma_a[\tau_1]=p^{-1}$ and $\E^\sigma_a[\tau_1^2]=(2-p)/p^2\leq 2/p^2$ for $p=\P^\sigma_a(\Delta_1\geq A_2)\geq C/\sqrt{A_2}$. Collecting the
above estimates and substituting them in (\ref{Hextest4}) then yields
$$
\E^\sigma_a[V(\sigma, L)] \leq C \sqrt{L} A_2^{\frac{7}{2}} = C L^{\frac{15}{16}},
$$
which by Markov inequality implies (\ref{Hextest2}) if $B_1$, and hence $L=e^{B_1\rho^{-\zeta}}$, is sufficiently large.
\qed

\section{Proof of Proposition \ref{P:coarsegrain}}\label{S:coarsegrain}

The deduction of Prop.~\ref{P:coarsegrain} from Prop.~\ref{P:cgZab} is model independent. Part of the proof is similar to its discrete time analogue (see e.g.\ the proofs of Proposition 2.3 and Lemma 2.4 in \cite{GLT09}), with the main difference being that in the integrals in (\ref{cgZab1})--(\ref{cgZab2}), we have excluded not
only contributions from $b\in [a, a+\eps L]$, but also from $b\in [L-\eps L, L]$. The latter requires new bounds.

First note that by (\ref{ctsz}), for any $0<a<b<L$ and $z\in [1, 1+L^{-1}]$, we have
\begin{eqnarray}
\E^Y_0[\Zi^{z, {\rm pin}}_{[a,b],Y}] &=& \sum_{m=1}^\infty \!\! \idotsint\limits_{ \atop \sigma_0=a<\sigma_1\cdots<\sigma_m=b} \!\!\!
z^m  \prod\limits_{i=1}^{m} K(\sigma_i-\sigma_{i-1}) d\sigma_1\cdots d\sigma_{m-1} \nn \\
&=& z (G^{X-Y})^{-1} \E^{X,Y}_{0,0}\big[e^{z (G^{X-Y})^{-1} \int_a^b 1_{\{X_s=Y_s\}}ds} 1_{\{X_b=Y_b\}} \big| X_a=Y_a \big] \nn \\
&\leq& z (G^{X-Y})^{-1} \E^{X,Y}_{0,0}\big[e^{(G^{X-Y})^{-1} (1+\int_a^b 1_{\{X_s=Y_s\}}ds)} 1_{\{X_b=Y_b\}} \big| X_a=Y_a \big] \nn \\
&=& z e^{(G^{X-Y})^{-1}} P(b-a) \leq C P(b-a), \label{cg1}
\end{eqnarray}
where $P(t)$ is defined in (\ref{Pt1}), and $C>0$ is uniform in $L\geq 1$ and $Y$'s jump rate $\rho \in [0,1]$.
For the inequality in the third line, we used $z\leq 1+L^{-1}$ and $b-a\leq L$, while the equality in the last line follows by
Taylor expanding as in (\ref{ctszfree}) and setting $z=1$, i.e.\ $\beta=1/G^{X-Y}$, in (\ref{ctsz}) and then averaging over $Y$.
Let $I=\{i_1<i_2<\cdots <i_k\}\subset \{1,\cdots, m\}$ for some $m\in\N$. Then by (\ref{ctsfIZINY}), Prop.~\ref{P:cgZab}, and (\ref{cg1}),
\begin{eqnarray}
&& \E^Y_0\big[f_I(Y) \Zi^{z,I}_{t,Y}\big]  \label{cg2}   \\
&\leq& \nn
\idotsint\limits_{a_j<b_j \in \Lambda_{i_j} \atop 1\leq j\leq k}
(Cz)^k \prod_{j=1}^k K(a_j-b_{j-1}) (C 1_{\{b_j-a_j<\eps L \}}+C1_{\{b_j\geq (i_j-\eps)L\}} +\delta)P(b_j-a_j) \prod_{j=1}^k {\rm d}a_j\, {\rm d}b_j.
\end{eqnarray}
We will show that for any $\eta>0$, if $\eps>0$ and $\delta>0$ are chosen sufficiently small, then we have
\be\label{cg2.5}
\E^Y_0\big[f_I(Y) \Zi^{z,I}_{t,Y}\big] \leq \idotsint\limits_{a_j<b_j \in \Lambda_{i_j} \atop 1\leq j\leq k}
(Cz\eta )^k \prod_{j=1}^k K(a_j-b_{j-1}) P(b_j-a_j) \prod_{j=1}^k {\rm d}a_j\, {\rm d}b_j.
\ee
Prop.~\ref{P:coarsegrain} then follows from the bound
\be\label{cg6}
P_L(I) := \idotsint\limits_{a_j<b_j \in \Lambda_{i_j} \atop 1\leq j\leq k} \prod_{j=1}^k K(a_j-b_{j-1}) P(b_j-a_j) \prod_{j=1}^k {\rm d}a_j\, {\rm d}b_j \leq C_L \prod_{j=1}^{|I|} \frac{C}{(i_j-i_{j-1})^{\frac{3}{2}}},
\ee
where $C_L$ depends only on $L$.

First we give a proof of (\ref{cg6}), which is similar to its discrete time counterpart, \cite[Lemma 2.4]{GLT09}. We partition $I$ into blocks of consecutive integers $\{u_1, u_1+1,\cdots, v_1\}$, $\{u_2, u_2+1,\cdots, v_2\}$, \ldots, $\{u_l,u_l+1,\cdots, v_l\}$, where $u_j-v_{j-1}\geq 2$ for all $2\leq j\leq l$. When substituting the definition of $P(b_j-a_j)$ in (\ref{Pt1})  into (\ref{cg6}), the resulting multifold expansion is the probability of a set of renewal configurations, where the variables of integration constitute the renewal configuration $\sigma$. By only retaining the constraint that $\sigma$ intersects $\Lambda_{u_i}$ and $\Lambda_{v_i}$ for $1\leq i\leq l$, we obtain
\be\label{cg7}
P_L(I) \leq \!\!\!\!\!\!\!\! \int\limits_{a_1<b_1 \atop a_1 \in \Lambda_{u_1}, b_1\in \Lambda_{v_1}}\!\!\!\!\!\cdots\!\!\!\!\! \int\limits_{a_l< b_l\atop a_l\in \Lambda_{u_l},b_l\in \Lambda_{v_l}} \!\!\!\!\! \prod_{j=1}^l K(a_j-b_{j-1}) P(b_j-a_j) \prod_{j=1}^l {\rm d}a_j\, {\rm d}b_j,
\ee
where $b_0:=0$. We integrate out one pair of variables $a_j, b_j$ at a time. For $j=l$,
\be\label{cg8}
\iint\limits_{a_l<b_l \atop a_l\in \Lambda_{u_l}, b_l\in \Lambda_{v_l}} K(a_l-b_{l-1})P(b_l-a_l) db_l da_l \leq \frac{C}{(u_l-v_{l-1})^{\frac{3}{2}}L^{\frac{3}{2}}}
\iint\limits_{a_l<b_l \atop a_l\in \Lambda_{u_l}, b_l\in \Lambda_{v_l}} \!\!\!\!\!\! P(b_l-a_l) db_l da_l,
\ee
since $a_l-b_{l-1}\geq (u_l-v_{l-1}-1)L \geq L$, and hence $K(a_l-b_{l-1})\leq \frac{C_1}{(u_l-v_{l-1}-1)^{\frac{3}{2}}L^{\frac{3}{2}}}\leq\frac{C}{(u_l-v_{l-1})^{\frac{3}{2}}L^{\frac{3}{2}}}$.
Regardless of whether $u_l=v_l$ or $v_l>u_l$, by Lemma \ref{L:renewalprob}, we have
\begin{eqnarray*}
\iint\limits_{a_l<b_l \atop a_l \in \Lambda_{u_l}, b_l\in \Lambda_{v_l}} \!\!\!\!\!\!\!\!  P(b_l-a_l) db_l da_l
\leq \!\!\!\!\!\!\!\! \iint\limits_{a_l<b_l \atop a_l \in \Lambda_{u_l}, b_l\in \Lambda_{v_l}} \!\!\!\!\!\!\!\! \frac{C}{\sqrt{b_l-a_l}} db_l da_l \leq C L^{\frac{3}{2}}.
\end{eqnarray*}
Integrating out $a_l, b_l$ in (\ref{cg7}) thus gives a factor $C(u_l-v_{l-1})^{-\frac{3}{2}}$. Iterating this procedure
then gives the bound in (\ref{cg6}), where a prefactor $C_L=L^{\frac{3}{2}}$ arises when we integrate out $a_1$ and $b_1$ in the case $u_1=1$. This proves (\ref{cg6}).

To deduce (\ref{cg2.5}) from (\ref{cg2}), we first bound the contributions from $C 1_{\{b_j-a_j<\eps L \}}$. We claim that there exists some $C>0$ depending only on
$K(\cdot)$ and uniform in $\rho\in [0,1]$, such that for all $L$ sufficiently large, $\eps\in (0,1/4)$, and $a\leq 0<L\leq b$, we have
\begin{eqnarray}
\int\limits_{0\leq s<t\leq L \atop t-s < \eps L} \!\!\!\!\!\! K(s-a) P(t-s)K(b-t) dtds &\leq& C\sqrt{\eps}\!\!\!\!\!\! \int\limits_{0\leq s<t\leq L}\!\!\!\!\!\! K(s-a) P(t-s)K(b-t) dtds, \label{cg3} \\
\int\limits_{0\leq s<t\leq L \atop t-s < \eps L} \!\!\!\!\!\! K(s-a) P(t-s) dtds &\leq& C\sqrt{\eps}\!\!\!\!\!\! \int\limits_{0\leq s<t\leq L}\!\!\!\!\!\! K(s-a) P(t-s) dtds.              \label{cg4}
\end{eqnarray}
To prove (\ref{cg3}), note that either $s\leq L/2$ or $s>L/2$ in the integral. Using the fact that $K(t)\sim \frac{C}{t^{\frac{3}{2}}}$ by the local central limit
theorem and the fact that $\int_0^t P(s)ds\sim C\sqrt{t}$ by Lemma \ref{L:renewalprob}, we have
$$
\int\limits_{0\leq s<t\leq L \atop s\leq \frac{L}{2}, t-s < \eps L} \!\!\!\!\!\!\!\!\! K(s-a) P(t-s)K(b-t) dtds
\leq \frac{C\sqrt{\eps L}}{(b-\frac{3L}{4})^{\frac{3}{2}}} \int\limits_{0\leq s\leq L/2} \!\!\!\!\!\! K(s-a)ds
\leq \frac{C\sqrt{\eps L}}{b^{\frac{3}{2}}} \int\limits_{0\leq s\leq L/2} \!\!\!\!\!\! K(s-a)ds,
$$
where we used $b-t>b-3L/4$ and $b\geq L$. On the other hand,
$$
\int\limits_{0\leq s<t\leq L} \!\!\!\!\!\! K(s-a) P(t-s)K(b-t) dtds \geq \!\!\!\!\!\! \int\limits_{0\leq s\leq \frac{L}{2} \atop 0<t-s\leq \frac{L}{4}}
\!\!\!\!\!\! K(s-a) P(t-s)K(b-t) dtds \geq \frac{C\sqrt{L}}{b^{\frac{3}{2}}} \!\!\!\!\! \int\limits_{0\leq s\leq L/2} \!\!\!\!\!\! K(s-a)ds.
$$
Together with a similar bound for the LHS of (\ref{cg3}) integrated over $s>L/2$, this implies (\ref{cg3}). The proof of (\ref{cg4}) is similar and will be
omitted. Substituting (\ref{cg3}) and (\ref{cg4}) into (\ref{cg2}) then gives

\be\label{cg5}
\E^Y_0\big[f_I(Y) \Zi^{z,I}_{t,Y}\big] \leq  (Cz )^k \idotsint\limits_{a_j<b_j \in \Lambda_{i_j} \atop 1\leq j\leq k}
\prod_{j=1}^k K(a_j-b_{j-1}) (C1_{\{b_j\geq (i_j-\eps)L \}} + \tilde\eta) P(b_j-a_j) \prod_{j=1}^k {\rm d}a_j\, {\rm d}b_j,
\ee
where $\tilde\eta = C\sqrt{\eps}+\delta$, which can be made arbitrarily small by choosing $\eps$ and $\delta$ small. By expanding
the product $\prod_{j=1}^k(C1_{\{b_j\geq (i_j-\eps)L \}} + \tilde\eta)$, we note that (\ref{cg2.5}) follows once we show that there exists some $C$
such that for any $J\subset \{1, \dots, k\}$, we have
\be\label{cg9}
\idotsint\limits_{1\leq j\leq k :\ a_j<b_j \in \Lambda_{i_j} \atop j\in J :\  b_j\geq (i_j-\eps)L}
\prod_{j=1}^k K(a_j-b_{j-1}) P(b_j-a_j) \prod_{j=1}^k {\rm d}a_j\, {\rm d}b_j \leq (C\sqrt{\eps})^{|J|} P_L(I),
\ee
where $P_L(I)$ was defined in (\ref{cg6}).

If $J=\emptyset$, then (\ref{cg9}) is trivial; otherwise, let $l$ be the largest element in $J$. It suffices to show that we can replace
the indicator $1_{\{b_l\geq (i_l-\eps)L\}}$ by the factor $C\sqrt{\eps}$. We can then apply the argument inductively to deduce (\ref{cg9}).
There are three cases: either (1) $l=k$; or (2) $i_{l+1}-i_l \geq 2$; or (3) $i_{l+1}-i_l=1$. For the case $l=k$, it suffices to show that uniformly
in $b_{k-1}\in \Lambda_{i_{k-1}}$, we have
\be\label{cg10}
\iint\limits_{a_k<b_k\in \Lambda_{i_k}  \atop b_k\geq (i_k-\eps)L} \!\!\!\!
K(a_k-b_{k-1}) P(b_k-a_k) da_k db_k \leq C\sqrt{\eps} \!\!\!\!\!
\iint\limits_{a_k<b_k\in \Lambda_{i_k}} \!\!\!\! K(a_k-b_{k-1}) P(b_k-a_k) da_k db_k.
\ee
Note that by Lemma \ref{L:renewalprob}, uniformly in $u>0$, we have $\int_u^{u+\eps L} P(s)ds \leq C\sqrt{\eps L}$ for $L$ large. Uniformly in $u>0$,
we also have $\int_u^{u+L} K(s)ds \leq 2 \int_u^{u+\frac{L}{2}} K(s)ds$. It follows that
$$
\iint\limits_{ a_k<b_k\in \Lambda_{i_k}  \atop b_k\geq (i_k-\eps)L} \!\!\!\!
K(a_k-b_{k-1}) P(b_k-a_k) da_k db_k \leq 2 C\sqrt{\eps L} \!\!\!\!\!\!\!\!\!\!\!\!\!\!\!
\int\limits_{(i_k-1)L \leq a_k\leq (i_k-\frac{1}{2})L}\!\!\!\!\!\!\!\!\!\!\!\!\!\!\! K(a_k-b_{k-1}) da_k.
$$
On the other hand, by Lemma \ref{L:renewalprob}, $\int_0^t P(s)ds \sim C\sqrt{t}$. Therefore for $L$ sufficiently large,
$$
\iint\limits_{ a_k<b_k\in \Lambda_{i_k} } \!\!\!\!\!\!
K(a_k-b_{k-1}) P(b_k-a_k) da_k db_k \geq \!\!\!\!\!\!\!\!\!\!\!\!\!\!\!\!\!\!\!
\int\limits_{(i_k-1)L \leq a_k\leq (i_k-\frac{1}{2})L \atop 0\leq b_k-a_k\leq \frac{L}{2}} \!\!\!\!\!\!\!\!\!\!\!\!\!\!\!\!\!\!\!
K(a_k-b_{k-1}) P(b_k-a_k) da_k db_k
\geq C\sqrt{L} \!\!\!\!\!\!\!\!\!\!\!\!\!\!\!\!\!\!\!\!\!
\int\limits_{(i_k-1)L \leq a_k\leq (i_k-\frac{1}{2})L}\!\!\!\!\!\!\!\!\!\!\!\!\!\!\!\!\!\!\!\!\!  K(a_k-b_{k-1}) da_k.
$$
The above two estimates together imply (\ref{cg10}).

For case (2), $i_{l+1}-i_l \geq 2$, it suffices to show that uniformly in $b_{l-1}\in \Lambda_{i_{l-1}}$ and $a_{l+1}\in \Lambda_{i_{l+1}}$, we have
$$
\iint\limits_{a_l<b_l\in \Lambda_{i_l}  \atop b_l\geq (i_l-\eps)L} \!\!\!\!
K(a_l-b_{l-1}) P(b_l-a_l) K(a_{l+1}-b_l) da_l db_l \leq C\sqrt{\eps} \!\!\!\!\!
\iint\limits_{a_l<b_l\in \Lambda_{i_l}} \!\!\!\! K(a_l-b_{l-1}) P(b_l-a_l) K(a_{l+1}-b_l)  da_l db_l.
$$
This follows from the same proof as for (\ref{cg10}) once we note that, because $a_{l+1}-b_l\geq L$, uniformly in $s_1,s_2\in \Lambda_{i_l}$ and
$t_1,t_2\in \Lambda_{i_{l+1}}$, we have $C\leq \frac{K(t_1-s_1)}{K(t_2-s_2)} \leq C^{-1}$ for some $C\in (0,\infty)$ depending only on $K(\cdot)$.

For case (3), $i_{l+1}-i_l=1$, there are two subcases: either $l+1=k$ or $l+1<k$. We only examine the case $l+1<k$, since the case $l+1=k$
is similar and simpler. To simplify notation, we will shift coordinates and assume $l=1$ and $i_l=1$. Since $l$ is the largest element in $J$,
it suffices to show that uniformly in $b_0\leq 0$ and $a_3\geq 2L$, we have
\begin{eqnarray}
\!\!\!\!\!\!\!\!\!\!\!\!\! && \idotsint\limits_{b_1\geq (1-\eps)L \atop 0<a_1<b_1<L<a_2<b_2<2L} \!\!\!\!\!\!\!\!\!\!\!\!\!
K(a_1-b_0) P(b_1-a_1) K(a_2-b_1) P(b_2-a_2) K(a_3-b_2) da_1 db_1 da_2 db_2  \nn \\
\!\!\!\!\!\!\!\!\!\!\!\!\! &\leq& C\sqrt{\eps} \!\!\!\!\!\!\!\!\!\!\!\!\! \idotsint\limits_{0<a_1<b_1<L<a_2<b_2<2L} \!\!\!\!\!\!\!\!\!\!\!\!\!
K(a_1-b_0) P(b_1-a_1) K(a_2-b_1) P(b_2-a_2) K(a_3-b_2) da_1 db_1 da_2 db_2. \label{cg11}
\end{eqnarray}
By restricting the region of integration to $a_1\in [0, L/4]$, $b_1\in [3L/4, L]$, $a_2\in [L, 5L/4]$ and $b_2\in[7L/4, 2L]$, and using the fact that
$P(t)\sim \frac{C}{\sqrt{t}}$, $K(t)\sim \frac{C}{t^{\frac{3}{2}}}$, $\int_t^\infty K(s)ds\sim \frac{C}{\sqrt t}$, we find
\begin{eqnarray}
&& \!\!\!\!\!\!\!\!\!\!\!\!\! \idotsint\limits_{0<a_1<b_1<L<a_2<b_2<2L} \!\!\!\!\!\!\!\!\!\!\!\!\!
K(a_1-b_0) P(b_1-a_1) K(a_2-b_1) P(b_2-a_2) K(a_3-b_2) da_1 db_1 da_2 db_2  \nn \\
&\geq& \frac{C}{\sqrt L} \int\limits_0^{\frac{L}{4}} K(a_1-b_0)da_1 \int\limits_{\frac{7L}{4}}^{2L} K(a_3-b_2)db_2. \label{cg12}
\end{eqnarray}
To upper bound the LHS of (\ref{cg11}), we claim that uniformly in all $b_1\leq L<2L\leq a_3$, we have
\be\label{cg13}
\int\limits_L^{2L} \int\limits_{a_2}^{2L}
K(a_2-b_1)P(b_2-a_2)K(a_3-b_2) db_2 da_2\leq \frac{C}{\sqrt L} \int\limits_L^{\frac{5L}{4}} K(a_2-b_1)da_2 \int\limits_{\frac{7L}{4}}^{2L} K(a_3-b_2)db_2,
\ee
and uniformly for all $b_0 \leq 0$ and $(1-\eps)L <b_1 <L$, we have
\be\label{cg14}
\int\limits_0^{b_1} K(a_1-b_0) P(b_1-a_1) da_1 \leq \frac{C}{\sqrt L} \int\limits_0^{\frac{L}{4}} K(a_1-b_0)da_1,
\ee
which when substituted into the LHS of (\ref{cg11}) imply that
\begin{eqnarray}
\!\!\!\!\!\!\!\!\!\!\!\!\! && \!\!\!\!\!\!\!\!\!\!\!\!\! \idotsint\limits_{b_1\geq (1-\eps)L \atop 0<a_1<b_1<L<a_2<b_2<2L} \!\!\!\!\!\!\!\!\!\!\!\!\!
K(a_1-b_0) P(b_1-a_1) K(a_2-b_1) P(b_2-a_2) K(a_3-b_2) da_1 db_1 da_2 db_2  \nn \\
\!\!\!\!\!\!\!\!\!\!\!\!\! &\leq& \frac{C}{L} \int\limits_{(1-\eps)L}^L  \int\limits_L^{\frac{5L}{4}}K(a_2-b_1) da_2 db_1
\int\limits_0^{\frac{L}{4}} K(a_1-b_0)da_1 \int\limits_{\frac{7L}{4}}^{2L} K(a_3-b_2)db_2 \nn \\
\!\!\!\!\!\!\!\!\!\!\!\!\! &\leq& \frac{C\sqrt{\eps}}{\sqrt{L}} \int\limits_0^{\frac{L}{4}} K(a_1-b_0)da_1 \int\limits_{\frac{7L}{4}}^{2L} K(a_3-b_2)db_2,
\end{eqnarray}
where we used the fact that $\int_L^{\frac{5L}{4}}K(a_2-b_1) da_2 \leq \frac{C}{\sqrt{L-b_1}}$. Together with (\ref{cg12}), this implies (\ref{cg11}).

To prove (\ref{cg13}), note that the bound therein certainly holds if we restrict integration to $a_2\in [L, \frac{5L}{4}]$ and $b_2\in [\frac{7L}{4}, 2L]$.
If either of the constraints on $a_2$ and $b_2$ fails, without loss of generality, say $a_2\in [\frac{5L}{4},2L]$, then because $K(t)\leq \frac{C}{t^{\frac{3}{2}}}$
and $\int_0^t P(s)ds\leq C\sqrt{t}$, we have
\begin{eqnarray*}
\!\!\!\!\!\! && \int\limits_{\frac{5L}{4}}^{2L} \int\limits_{a_2}^{2L} K(a_2-b_1)P(b_2-a_2)K(a_3-b_2) db_2 da_2 \leq \frac{C\sqrt{L}}{(\frac{5L}{4}-b_1)^{\frac{3}{2}}}  \int_L^{2L} K(a_3-b_2)db_2 \\
&\leq&  \frac{C\sqrt{L}}{(\frac{5L}{4}-b_1)^{\frac{3}{2}}} \int_{\frac{7L}{4}}^{2L} K(a_3-b_2)db_2
\leq \frac{C}{\sqrt L} \int\limits_L^{\frac{5L}{4}} K(a_2-b_1)da_2 \int\limits_{\frac{7L}{4}}^{2L} K(a_3-b_2)db_2,
\end{eqnarray*}
since $\int_L^{\frac{5L}{4}} K(a_2-b_1)da_2 \geq \frac{C L}{4} (\frac{5L}{4}-b_1)^{-\frac{3}{2}}$. This proves (\ref{cg13}). The proof of (\ref{cg14}) is
similar and will be omitted. This completes the proof of (\ref{cg9}) as well as of Prop.~\ref{P:coarsegrain}.
\qed

\section{Proof of Lemma \ref{L:HLYbound}}\label{S:HLYbound}
Note that (\ref{ctsEHN}) is obvious. For $s\in [0,L]$, let us denote
$$
h_L(s,Y) := \int\limits_{s<t<L \atop A_1<t-s<A_2} \frac{1_{\{Y_s=Y_t\}}}{(\log (t-s))^\xi} dt.
$$
Then
\begin{eqnarray*}
{\rm Var}(H_L(Y)) &=& 2 \!\!\!\!\!\!\!\! \int\limits_{0<s_1<s_2<L}\!\!\!\!\!\!\!\!\!\! \big(\E^Y_0[h_L(s_1,Y) h_L(s_2,Y)] - \E^Y_0[h_L(s_1,Y)]\E^Y_0[h_L(s_2,Y)]\big) ds_1 ds_2 \\
&\leq& 2 \!\!\!\!\!\!\!\! \iiiint\limits_{0<s_1<s_2<L, \, t_1,t_2<L \atop A_1<t_1-s_1, t_2-s_2<A_2} \!\!\!\!\!\!\!\!\!\!
\frac{\big| \P^Y_0(Y_{s_1}=Y_{t_1}, Y_{s_2}=Y_{t_2})-\P^Y_0(Y_{s_1}=Y_{t_1})\P^Y_0(Y_{s_2}=Y_{t_2})\big|}{(\log (t_1-s_1) \log (t_2-s_2))^\xi} dt_1dt_2ds_1ds_2 \\
&\leq & 2 \int\limits_{0<s_1<s_2<L} \phi(s_2-s_1) ds_1 ds_2 \leq 2 L \int_0^\infty \phi(w) dw,
\end{eqnarray*}
where
\begin{eqnarray}
\phi(w) &=& \int_{A_1}^\infty \int_{w+A_1}^\infty \frac{\big| \P^Y_0(Y_0=Y_{s_1}, Y_w=Y_{s_2})-\P^Y_0(Y_0=Y_{s_1})\P^Y_0(Y_w=Y_{s_2})\big|}{(\log s_1 \log (s_2-w))^\xi} ds_2 ds_1.
 \label{phidec}
\end{eqnarray}
To prove (\ref{ctsvarHN}), it suffices to show that $\int_0^\infty \phi(w)dw \leq C/\rho^3$.

Note that in (\ref{phidec}), $s_1,s_2$ fall into three regions: (0) $0< s_1 <w$; (1) $w<s_1<s_2$; (2) $w<s_2<s_1$. In case (0), the integrand in (\ref{phidec}) is
$0$ by the independent increment properties of $Y$. In case (1), let $r_1=s_1-w$ and $r_2=s_2-s_1$, while in case (2) let $r_1=s_2-w$ and $r_2=s_1-s_2$, then
\be
\phi(w) = {\rm I}(w) +{\rm II}(w)
\ee
with
$$
\begin{aligned}
{\rm I}(w) & = \iint\limits_{[0,\infty)^2}  1_{\{w+r_1>A_1, r_1+r_2>A_1\}} \frac{\big|\P^Y_0(Y_0=Y_{w+r_1}, Y_w=Y_{w+r_1+r_2}) -p_{\rho(w+r_1)}(0)p_{\rho(r_1+r_2)}(0)\big|}{(\log(w+r_1)\log(r_1+r_2))^\xi}\, dr_1 dr_2, \\
{\rm II}(w) & = \iint\limits_{[0,\infty)^2} 1_{\{r_1>A_1\}}
\frac{\big|\P^Y_0(Y_0=Y_{w+r_1+r_2}, Y_w=Y_{w+r_1}) -p_{\rho(w+r_1+r_2)}(0)p_{\rho r_1}(0)\big|}{(\log(w+r_1+r_2)\log r_1)^\xi}\, dr_1 dr_2.
\end{aligned}
$$
Since $\xi>1/2$, we establish (\ref{ctsvarHN}) once we show that there exists $C>0$ such that
\be\label{IIIwbound}
\begin{aligned}
{\rm I}(w), {\rm II}(w)\ &\leq \ \ \ \ \ \  \frac{C}{\rho^2}             \qquad & &\mbox{for all } w>0, \\
{\rm I}(w), {\rm II}(w)\ &\leq \frac{C}{\rho^3 w(\log w)^{2\xi}} \qquad & &\mbox{for all } w>A_1=e.
\end{aligned}
\ee

In ${\rm I}(w)$, by Lemmas \ref{L:LCLT} and \ref{L:est1},
\begin{eqnarray*}
&& \P^Y_0(Y_0=Y_{w+r_1}, Y_w=Y_{w+r_1+r_2}) \\
&=& \sum_{x\in\Z^3} p_{\rho w}(x) p_{\rho r_1}(-x) p_{\rho r_2}(x) \leq
\min\Big\{ p_{\rho r_1}(0)p_{\rho r_2}(0), \frac{C}{\rho^3 (wr_1+wr_2+r_1r_2)^{\frac{3}{2}}}\Big\},
\end{eqnarray*}
from which we easily deduce that ${\rm I}(w)\leq 2(\int_0^\infty p_{\rho r}(0)dr)^2 = 2G^2 \rho^{-2}$.
Similarly, ${\rm II}(w) \leq 2 G^2 \rho^{-2}$. On the other hand, by the local central limit theorem, Lemma \ref{L:LCLT},
we have
\be\label{Iw}
{\rm I}(w) \leq \frac{C}{\rho^3} \iint_{[0,\infty)^2} 1_{\{w+r_1>A_1, r_1+r_2>A_1\}} \frac{\frac{1}{(wr_1 + wr_2 +r_1 r_2)^{\frac{3}{2}}}+\frac{1}{(w+r_1)^{\frac{3}{2}} (r_1+r_2)^{\frac{3}{2}}} }{(\log(w+r_1)\log(r_1+r_2))^\xi} \, dr_1 dr_2.
\ee
Let $r_1 = w t_1$ and $r_2 = w t_2$, and assume $w>A_1$, then (\ref{Iw}) becomes
\begin{eqnarray}
&& {\rm I}(w) \leq \frac{C}{\rho^3 w} \iint_{[0,\infty)^2} 1_{\{1+t_1>A_1 w^{-1}, t_1+t_2>A_1w^{-1}\}}
\frac{\frac{1}{(t_1 + t_2 +t_1 t_2)^{\frac{3}{2}}}+ \frac{1}{(1+t_1)^{\frac{3}{2}} (t_1+t_2)^{\frac{3}{2}}} }{(\log(w(1+t_1))\log(w(t_1+t_2)))^\xi} \, dt_1 dt_2  \nonumber \\
&\leq&
\frac{C}{\rho^3 w (\log w)^{2\xi}}      \label{Iw1}
\!\! \iint\limits_{t_1, t_2\geq 0 \atop t_1+t_2\geq 1} \!\!\! \Big(\frac{1}{(t_1+t_2+t_1t_2)^{\frac{3}{2}}}+\frac{1}{(1+t_1)^{\frac{3}{2}}(t_1+t_2)^{\frac{3}{2}}}\Big) dt_1dt_2 \\
&& \ +\ \frac{C}{\rho^3 w(\log w)^\xi}  \label{Iw2}
\!\!\! \iint\limits_{t_1, t_2\geq 0 \atop t_1+t_2\leq 1} 1_{\{t_1+t_2>A_1w^{-1}\}}\frac{2}{(t_1+t_2)^{\frac{3}{2}}(\log(w(t_1+t_2)))^\xi} dt_1 dt_2.
\end{eqnarray}
The integral in (\ref{Iw1}) is finite. Letting $y=t_1+t_2$, the integral in (\ref{Iw2}) equals
\begin{eqnarray}
\int_{\frac{A_1}{w}}^1 \frac{2 }{\sqrt{y} (\log(wy))^\xi} dy &=& \frac{1}{\sqrt w}\int_{A_1}^w \frac{2}{\sqrt{x} (\log x)^\xi} dx  \label{Iw3}\\
&=&
\frac{4\sqrt{x}}{\sqrt{w}(\log x)^\xi}\Big|_{A_1}^w + \frac{4\xi}{\sqrt w}\int_{A_1}^w \frac{1}{\sqrt{x} (\log x)^{1+\xi}}dx \leq \frac{C}{(\log w)^\xi}, \nonumber
\end{eqnarray}
which proves that ${\rm I}(w) \leq \frac{C}{\rho^3 w(\log w)^{2\xi}}$ for $w>A_1$.

In ${\rm II}(w)$,
$$
\P^Y_0(Y_0=Y_{w+r_1+r_2}, Y_w=Y_{w+r_1}) = \sum_{x\in\Z^3} p_{\rho w}(x) p_{\rho r_1}(0) p_{\rho r_2}(-x) = p_{\rho r_1}(0)p_{\rho(w+r_2)}(0).
$$
Therefore
\be
{\rm II}(w) = \iint\limits_{[0,\infty)^2} 1_{\{r_1>A_1\}}
\frac{p_{\rho r_1}(0) |p_{\rho(w+r_2)}(0)-p_{\rho(w+r_1+r_2)}(0)|} {(\log(w+r_1+r_2)\log r_1)^\xi}\, dr_1 dr_2. \label{IIw1}
\ee
We separate the integral in (\ref{IIw1}) according to whether $r_1>w$ or $r_1<w$. When $r_1>w$, we have
\begin{eqnarray*}
&& \int_w^\infty\int_0^\infty 1_{\{r_1>A_1\}}
\frac{p_{\rho r_1}(0)|p_{\rho(w+r_2)}(0)-p_{\rho(w+r_1+r_2)}(0)|} {(\log(w+r_1+r_2)\log r_1)^\xi} dr_2 dr_1  \\
&\leq&
\frac{C\sqrt{w}}{\rho^{\frac{3}{2}}(\log w)^\xi} \int_1^\infty \int_0^\infty  \frac{p_{\rho w(1+t_2)}(0)+p_{\rho w(1+t_1+t_2)}(0)} {t_1^{\frac{3}{2}}(\log (wt_1))^\xi}\, dt_2 dt_1 \\
&\leq& \frac{C\sqrt{w}}{\rho^3 (\log w)^\xi} \int_1^\infty\int_0^\infty \frac{\frac{1}{(w(1+t_2))^{\frac{3}{2}}}}{t_1^{\frac{3}{2}}(\log w)^\xi} dt_2 dt_1 \\
&\leq& \frac{C}{\rho^3 w(\log w)^{2\xi}},
\end{eqnarray*}
where we used the local central limit theorem to bound $p_s(0)\leq C s^{-\frac{3}{2}}$ and made the change of variables $r_1=wt_1$ and $r_2=wt_2$.
When $0<r_1<w$ in (\ref{IIw1}), by Lemma \ref{L:est2}, we have
\begin{eqnarray*}
&& \int_0^w\int_0^\infty 1_{\{r_1>A_1\}}
\frac{p_{\rho r_1}(0)|p_{\rho(w+r_2)}(0)-p_{\rho(w+r_1+r_2)}(0)|} {(\log(w+r_1+r_2)\log r_1)^\xi} dr_2 dr_1  \\
&\leq& \frac{C}{(\log w)^\xi} \int_0^w\int_0^\infty 1_{\{r_1>A_1\}} \frac{p_{\rho r_1}(0)\frac{r_1}{\rho^{\frac{3}{2}}(w+r_2)^{\frac{5}{2}}}}{(\log r_1)^\xi} dr_2 dr_1 \\
&\leq& \frac{C}{\rho^3 (\log w)^\xi}  \int_0^w\int_0^\infty 1_{\{r_1>A_1\}} \frac{1}{\sqrt{r_1}(w+r_2)^{\frac{5}{2}}(\log r_1)^\xi }dr_2 dr_1 \\
&=& \frac{C}{\rho^3 w(\log w)^\xi}  \int_0^1\int_0^\infty 1_{\{t_1>A_1w^{-1}\}} \frac{1}{\sqrt{t_1}(1+t_2)^{\frac{5}{2}} (\log(wt_1))^\xi}dt_2 dt_1 \\
&\leq& \frac{C}{\rho^3 w(\log w)^\xi}  \int_{\frac{A_1}{w}}^1 \frac{1}{\sqrt{t_1}(\log(wt_1))^\xi}dt_1 \leq \frac{C}{\rho^3 w(\log w)^{2\xi}},
\end{eqnarray*}
where the last inequality follows from the same calculation as in (\ref{Iw3}). This proves that ${\rm II}(w) \leq \frac{C}{\rho^3 w(\log w)^{2\xi}}$ for $w>A_1$ and concludes the proof of Lemma \ref{L:HLYbound}.
\qed

\section{Proof of Lemma \ref{L:Hsigmabd}}\label{S:Hsigmabd}

\noindent
{\bf Proof of Lemma \ref{L:Hsigmabd}~(\ref{Hsigmext})--(\ref{Hsigmext2}).} By definition, conditioned on $\sigma$,
$$
\E^{Y^\sigma}_0[H_{[\sigma_i, \sigma_{i+1}]}^{\rm ext}(Y^\sigma)] - \E^Y_0[H_{[\sigma_i, \sigma_{i+1}]}^{\rm ext}(Y)] = \!\!\!\!\!\!\!\! \iint\limits_{\sigma_i<s_1<\sigma_{i+1}<s_2 \atop A_1<s_2-s_1<A_2} \!\!\!\!\!\!
\frac{\P(Y^\sigma_{s_1}=Y^\sigma_{s_2}) - \P(Y_{s_1}=Y_{s_2})}{\log(s_2-s_1)} ds_2 ds_1.
$$
To prove (\ref{Hsigmext}), it suffices to show that
$$
\P(Y^\sigma_{s_1}=Y^\sigma_{s_2}) > \P(Y_{s_1}=Y_{s_2}).
$$
This follows from Lemma \ref{L:retprocomp}. Indeed, we can decompose $Y^\sigma_{s_2}-Y^\sigma_{s_1}$ as the
sum of independent random variables $Z_1, Z_2,\cdots, Z_{n+1}$, where $n$ is such that $\sigma_{i+n}\leq s_2<\sigma_{i+n+1}$, $Z_1=Y^\sigma_{\sigma_{i+1}}-Y^\sigma_{s_1}$, $Z_j=Y^\sigma_{\sigma_{i+j}}-Y^\sigma_{\sigma_{i+j-1}}$ for $2\leq i\leq n$, and $Z_{n+1}=Y^\sigma_{s_2}-Y^\sigma_{\sigma_{i+n}}$. From the definition of $Y^\sigma$,
$$
\P(Z_1=y) =\frac{\sum_x p^Y_{s_1-\sigma_i}(x)p^Y_{\sigma_{i+1}-s_1}(y)p^X_{\sigma_{i+1}-\sigma_i}(x+y)}{p^{X-Y}_{\sigma_{i+1}-\sigma_i}(0)} = \frac{p_{\rho(\sigma_{i+1}-s_1)}(y)p_{\sigma_{i+1}-\sigma_i+\rho(s_1-\sigma_{i})}(y)}{p_{(1+\rho)(\sigma_{i+1}-\sigma_i)}(0)},
$$
where we used the fact that $X$ and $Y$ have the same symmetric jump probability kernel with respective rates $1$ and $\rho$. Therefore $Z_1$ is distributed as $X_{\rho(\sigma_{i+1}-s_1)}$ conditioned on $X_{(1+\rho)(\sigma_{i+1}-\sigma_i)}=0$. Similarly, $Z_j$ for $2\leq j\leq n$ is distributed as $X_{\rho(\sigma_{i+j}-\sigma_{i+j-1})}$ conditioned on $X_{(1+\rho)(\sigma_{i+j}-\sigma_{i+j-1})}=0$, and $Z_{n+1}$ is distributed as $X_{\rho(s_2-\sigma_{i+n})}$
conditioned on $X_{(1+\rho)(\sigma_{i+n+1}-\sigma_{i+n})}=0$. Therefore Lemma \ref{L:retprocomp} applies. The proof of (\ref{Hsigmext2}) is analogous and simpler.
\qed
\bigskip

\noindent
{\bf Proof of Lemma \ref{L:Hsigmabd}~(\ref{Hsigmint}).} Without loss of generality, assume that $\sigma_0=a=0$, and let
$\sigma_{1}=\Delta$. For $0\leq s_1\leq s_2\leq \Delta$, we have
$$
\P(Y^\sigma_{s_1} = Y^\sigma_{s_2}) = \frac{\sum_{x,y\in\Z^3} p_{\rho s_1}(x) p_{\rho(s_2-s_1)}(0)p_{\rho(\Delta-s_2)}(y)p_{\Delta}(x+y)}{p_{(1+\rho)\Delta}(0)}
= \frac{p_{(1+\rho)\Delta-\rho(s_2-s_1)}(0)p_{\rho(s_2-s_1)}(0)}{p_{(1+\rho)\Delta}(0)}.
$$
Therefore, conditioned on $\sigma_0=0$ and $\sigma_{1}=\Delta$,
\begin{eqnarray}
&& \E^{Y^\sigma}_0[H_{[\sigma_0, \sigma_{1}]}^{\rm int}(Y^\sigma)] - \E^Y_0[H_{[\sigma_0, \sigma_{1}]}^{\rm int}(Y)]
= \!\!\!\! \iint\limits_{0<s_1<s_2<\Delta \atop A_1<s_2-s_1<A_2} \!\!\!\! \frac{\P( Y^\sigma_{s_1}=Y^\sigma_{s_2})-\P(Y_{s_1}=Y_{s_2})}{(\log(s_2-s_1))^\xi} ds_2 ds_1 \nonumber \\
&=& \!\!\!\! \iint\limits_{0<s_1<s_2<\Delta \atop A_1<s_2-s_1<A_2} \!\!\!\!
\frac{p_{\rho(s_2-s_1)}(0)(p_{(1+\rho)\Delta-\rho(s_2-s_1)}(0)-p_{(1+\rho)\Delta}(0))}{p_{(1+\rho)\Delta}(0) (\log(s_2-s_1))^\xi} ds_2 ds_1 \nonumber \\
&\geq& C \iint\limits_{0<s_1<s_2<\Delta \atop A_1<s_2-s_1<A_2}
\frac{((1+\rho)\Delta)^{\frac{3}{2}}\frac{\rho(s_2-s_1)}{((1+\rho)\Delta)^{\frac{5}{2}}}}{\rho^{\frac{3}{2}}(s_2-s_1)^{\frac{3}{2}} (\log(s_2-s_1))^\xi} ds_2 ds_1
\nonumber \\
&\geq&  \frac{C}{\Delta\sqrt \rho}  \iint\limits_{0<s_1<s_2<\Delta \atop A_1<s_2-s_1<A_2} \!\!\!\!\!
\frac{ds_2 ds_1}{ \sqrt{s_2-s_1}(\log(s_2-s_1))^\xi}
\geq  \frac{C\sqrt{\Delta}}{8\sqrt{\rho} (\log \Delta)^\xi} 1_{\{2A_1<\Delta<A_2\}},  \label{Hdiff}
\end{eqnarray}
where we have applied Lemma \ref{L:est2} and used the local central limit theorem to bound $p_t(0)\leq C_1t^{-\frac{3}{2}}$ for all $t\geq 0$ and $p_t\geq C_2t^{-\frac{3}{2}}$ for all $t\geq 1$. This proves (\ref{Hsigmint}).
\qed
\bigskip

\noindent
{\bf Proof of Lemma \ref{L:Hsigmabd}~(\ref{Hsigmavar}).} Without loss of generality, let $\sigma_0=a=0$ and $\sigma_{1}=\Delta$. We have
\begin{eqnarray}
&& {\rm Var}(H_{[\sigma_0, \sigma_{1}]}^{\rm int}(Y^\sigma)|\sigma) \nn \\
&=&  \!\!\!\!\!\!
\iint\limits_{0<s_1<s_2<\Delta, A_1<s_2-s_1<A_2, \atop 0<s'_1<s'_2<\Delta, A_1<s'_2-s'_1<A_2} \!\!\!\!\!\!\!\!\!\!\!\!\!\!\!\!\!\!
\frac{\P(Y^\sigma_{s_1}=Y^\sigma_{s_2}, Y^\sigma_{s'_1}=Y^\sigma_{s'_2})-\P(Y^\sigma_{s_1}=Y^\sigma_{s_2})\P(Y^\sigma_{s'_1}=Y^\sigma_{s'_2})}{(\log(s_2-s_1) \log(s'_2-s'_1))^\xi} ds_1 ds_2 ds'_1 ds'_2 \nonumber \\
&\leq& 2 \!\!\!\!\!\!\!\!
\iint\limits_{0<s_1<s_2<\Delta, A_1<s_2-s_1, \atop s_1<s'_1<s'_2<\Delta, A_1<s'_2-s'_1} \!\!\!\!\!\!\!\!\!\!\!\!\!\!\!
\frac{\big|\P(Y^\sigma_{s_1}=Y^\sigma_{s_2}, Y^\sigma_{s'_1}=Y^\sigma_{s'_2})-\P(Y^\sigma_{s_1}=Y^\sigma_{s_2})\P(Y^\sigma_{s'_1}=Y^\sigma_{s'_2})\big|}{(\log(s_2-s_1) \log(s'_2-s'_1))^\xi} ds_1 ds_2 ds'_1 ds'_2. \label{Htvb1}
\end{eqnarray}
In the integral above, $s_1, s_2, s'_1$ and $s'_2$ fall into three regions: (1) $s_1 <s_2<s'_1<s'_2$; (2) $s_1<s'_1<s_2<s'_2$; (3) $s_1<s'_1<s'_2<s_2$.
In region (1), let $r_1=s_2-s_1$, $r_2=s'_1-s_2$, $r_3=s'_2-s'_1$, and similarly in regions (2) and (3), let $r_1, r_2$ and $r_3$ be the successive increments of
the ordered variables. Let (1), (2) and (3) also denote the respective contributions to the integral in (\ref{Htvb1}) from the three regions. Then for (1),
we have
\begin{eqnarray*}
\P(Y^\sigma_{s_1}=Y^\sigma_{s_2}, Y^\sigma_{s'_1} = Y^\sigma_{s'_2}) \!\!\! &=&\!\!\!\! \frac{\sum_{x,y,z\in\Z^d} p_{\rho s_1}(x) p_{\rho r_1}(0) p_{\rho r_2}(y)p_{\rho r_3}(0)p_{\rho(\Delta-s_1-r_1-r_2-r_3)}(z)p_{\Delta}(x+y+z)}{p_{(1+\rho)\Delta}(0)} \\
\!\!\! &=& \!\!\!\! \frac{p_{\rho r_1}(0)p_{\rho r_3}(0)p_{(1+\rho)\Delta-\rho(r_1+r_3)}(0)}{p_{(1+\rho)\Delta}(0)},
\end{eqnarray*}
and
$$
\P(Y^\sigma_{s_1}=Y^\sigma_{s_2})\P(Y^\sigma_{s'_1}=Y^\sigma_{s'_2}) = \frac{p_{\rho r_1}(0)p_{(1+\rho)\Delta-\rho r_1}(0)p_{\rho r_3}(0)p_{(1+\rho)\Delta- \rho r_3}(0)}{p_{(1+\rho)\Delta}(0)^2}.
$$
Therefore
\begin{eqnarray}
&&\!\!\!\!\!\!\!\!\! (1) = \label{Htvb(1)} \\
&& \!\!\!\!\!\!\!\!\!\!
\iint_{0<s_1,r_1, r_2, r_3<\Delta,\ A_1<r_1, r_3 \atop s_1+r_1+r_2+r_3<\Delta}
\!\!\!\!\!\!\!\!\!\!\!\!\!\!\!\!\!\!\!\!\!\!\!
\frac{\big|\frac{p_{\rho r_1}(0)p_{\rho r_3}(0)p_{(1+\rho)\Delta-\rho(r_1+r_3)}(0)}{p_{(1+\rho)\Delta}(0)}-\frac{p_{\rho r_1}(0)p_{(1+\rho)\Delta-\rho r_1}(0)p_{\rho r_3}(0)p_{(1+\rho)\Delta- \rho r_3}(0)}{p_{(1+\rho)\Delta}(0)^2}\big|}
{(\log r_1 \log r_3)^\xi} ds_1 dr_1 dr_2 dr_3.  \nn
\end{eqnarray}
By similar considerations, we have
\begin{eqnarray}
\!\!\!\!\!\! &&\!\!\!\!\!\!\!\!\!\!\!\!\! (2) =   \label{Htvb(2)} \\
&&\!\!\!\!\!\!\!\!\!\!\!\!\!\!\!\!\!\!\!\!
\iint_{0<s_1,r_1, r_2, r_3<\Delta, \atop A_1<r_1+r_2, r_2+r_3, s_1+r_1+r_2+r_3<\Delta}
\!\!\!\!\!\!\!\!\!\!\!\!\!\!\! \frac{ds_1 dr_1 dr_2 dr_3}{(\log (r_1+r_2) \log (r_2+r_3))^\xi}
\Big| \frac{\sum_{x\in\Z^3}p_{\rho r_1}(x)p_{\rho r_2}(x)p_{\rho r_3}(x)p_{(1+\rho)\Delta-\rho(r_1+r_2+r_3)}(x)}{p_{(1+\rho)\Delta}(0)} \nn \\
&& \qquad \qquad \qquad \qquad \qquad
-\quad  \frac{p_{\rho(r_1+r_2)}(0)p_{(1+\rho)\Delta-\rho(r_1+r_2)}(0)p_{\rho(r_2+r_3)}(0)p_{(1+\rho)\Delta-\rho(r_2+r_3)}(0)}{p_{(1+\rho)\Delta}(0)^2}\Big|,   \nonumber\\
\!\!\!\!\!\! &&\!\!\!\!\!\!\!\!\!\!\!\!\! (3) = \label{Htvb(3)} \\
&&\!\!\!\!\!\!\!\!\!\!\!\!\!\!\!\!\!\!
\iint_{0<s_1,r_1, r_2, r_3<\Delta,\atop A_1<r_2, \, s_1+r_1+r_2+r_3<\Delta}
 \frac{ds_1 dr_1 dr_2 dr_3}{(\log (r_1+r_2+r_3) \log r_2)^\xi}
\Big| \frac{p_{\rho(r_1+r_3)}(0)p_{\rho r_2}(0)p_{(1+\rho)\Delta-\rho(r_1+r_2+r_3)}(0)}{p_{(1+\rho)\Delta}(0)} \nn \\
&& \qquad \qquad \qquad \qquad \qquad - \quad \frac{p_{\rho(r_1+r_2+r_3)}(0)p_{(1+\rho)\Delta- \rho(r_1+r_2+r_3)}(0)p_{\rho r_2}(0)p_{(1+\rho)\Delta-\rho r_2}(0)}{p_{(1+\rho)\Delta}(0)^2}\Big|.  \nonumber
\end{eqnarray}
We will show that (1), (2), (3) are all bounded by $C\rho^{-3}\Delta$ for some $C$ uniform in $\rho\in (0,1]$ and $\Delta>0$.

For (1), we have
\begin{eqnarray*}
(1) &\leq& \!\! \Delta^2 \!\!\!\!\!\! \iint\limits_{A_1<r_1, r_3<\Delta \atop r_1+r_3<\Delta} \!\!\!\!
\frac{p_{\rho r_1}(0)p_{\rho r_3}(0)\big|p_{(1+\rho)\Delta}(0)p_{(1+\rho)\Delta-\rho(r_1+r_3)}(0)- p_{(1+\rho)\Delta-\rho r_1}(0)p_{(1+\rho)\Delta-\rho r_3}(0)\big|}
{p_{(1+\rho)\Delta}(0)^2(\log r_1 \log r_3)^\xi} dr_1 dr_3 \\
&\leq& \frac{C \Delta^5}{\rho^3} \iint\limits_{A_1<r_1, r_3<\Delta}
\frac{\frac{\rho^2 r_1r_3}{\Delta^5}}
{r_1^\frac{3}{2}r_3^{\frac{3}{2}}(\log r_1 \log r_3)^\xi} dr_1 dr_3  \\
&\leq& \frac{C \Delta}{\rho},
\end{eqnarray*}
where we used the local central limit theorem to bound $p_s(0)$ and we applied Lemma \ref{L:est3} with $t=(1+\rho)\Delta- \rho(r_1+r_3)$.

We can bound (2) by passing the absolute value in (\ref{Htvb(2)}) inside. By Lemma \ref{L:LCLT}, $p_s(0)\leq \frac{C}{(1+s)^{\frac{3}{2}}}$ for some $C>0$ for all $s\geq 0$, and for $\Delta$ large, we have
\be\label{pDeltacomp}
\frac{p_{\Delta+s}(x)}{p_{(1+\rho)\Delta}(0)} < (1+\rho)^2 \qquad \mbox{for all } \ 0\leq s\leq \rho \Delta \ \ \mbox{and} \ \ x\in\Z^3.
\ee
Therefore
\begin{eqnarray}
&&  \iint\limits_{0<s_1,r_1, r_2, r_3<\Delta,\atop A_1<r_1+r_2, r_2+r_3, \, s_1+r_1+r_2+r_3<\Delta} \!\!\!\!\!\!\!\!
\frac{\frac{p_{\rho(r_1+r_2)}(0)p_{(1+\rho)\Delta- \rho(r_1+r_2)}(0)p_{\rho(r_2+r_3)}(0)p_{(1+\rho)\Delta-\rho(r_2+r_3)}(0)}{p_{(1+\rho)\Delta}(0)^2}}{(\log (r_1+r_2) \log (r_2+r_3))^\xi} ds_1 dr_1 dr_2 dr_3 \nonumber\\
&\leq& \frac{C \Delta}{\rho^3} \iint\limits_{0\leq r_1, r_2, r_3<\infty} \frac{1}{(1+r_1+r_2)^{\frac{3}{2}}(1+r_2+r_3)^{\frac{3}{2}}(\log (e+r_1+r_2) \log (e+r_2+r_3))^\xi} dr_1 dr_2 dr_3 \nonumber \\
&\leq& \frac{C\Delta}{\rho^3}, \label{Htvb2}
\end{eqnarray}
where the last integral is finite since integrating out $r_1$ and $r_3$ leads to a bound of the form $\int_0^\infty \frac{C}{(1+r_2)(\log (e+r_2))^{2\xi}}dr_2<\infty$. For the remaining term in (\ref{Htvb(2)}), we have
\begin{eqnarray}
&&  \iint\limits_{0<s_1,r_1, r_2, r_3<\Delta,\atop A_1<r_1+r_2, r_2+r_3, \, s_1+r_1+r_2+r_3<\Delta}  \!\!\!\!\!\!\!\!\!\!\!\!\!\!
\frac{\sum_{x\in\Z^3}p_{\rho r_1}(x)p_{\rho r_2}(x)p_{\rho r_3}(x)p_{(1+\rho)\Delta-\rho(r_1+r_2+r_3)}(x)}{p_{(1+\rho)\Delta}(0)(\log (r_1+r_2) \log (r_2+r_3))^\xi} ds_1 dr_1 dr_2 dr_3 \nonumber \\
&\leq& \frac{C\Delta}{\rho^3} \iint\limits_{0\leq r_1, r_2, r_3<\infty} \frac{dr_1 dr_2 dr_3}{(1+r_1r_2 +r_1r_3+r_2r_3)^{\frac{3}{2}}(\log (e+r_1+r_2) \log (e+r_2+r_3))^\xi} , \nonumber \\
&\leq& \frac{C\Delta}{\rho^3} \iint\limits_{0\leq r_1, r_2<\infty} \frac{dr_1 dr_2}{(\log(e+r_2))^{2\xi} (r_1+r_2) \sqrt{1+r_1r_2} }  \nonumber \\
&=& \frac{C \Delta}{\rho^3} \iint\limits_{0\leq t, r_2<\infty} \frac{dt dr_2}{(\log(e+r_2))^{2\xi} (1+t) \sqrt{1+t r_2^2} }  , \label{Htvb3}
\end{eqnarray}
where we used (\ref{pDeltacomp}), applied Lemma \ref{L:est1}, and made a change of variable $r_1=t r_2$. The integral in (\ref{Htvb3}) is clearly finite when integrated over $r_2>1$, since we can bound $\frac{1}{\sqrt{1+t r_2^2}}$ by $\frac{1}{r_2\sqrt{t}}$. For $0<r_2<1$, note that
$$
\int_0^\infty \frac{dt}{(1+t)\sqrt{1+tr_2^2}} = \int_0^\infty \frac{dw}{(w+r_2^2)\sqrt{1+w}}  \leq C -2 \ln r_2,
$$
which is integrable over $r_2\in [0,1]$. Therefore the integral in (\ref{Htvb3}) is finite, and together with (\ref{Htvb2}), this  shows that $(2)\leq C\rho^{-3}\Delta$.

For (3), we have
\begin{eqnarray*}
(3) &\leq&  \Delta \!\!\!\!\!\!\!\!\!\!\!\!\!\!
\iint\limits_{0<r_1, r_3<\Delta \atop A_1<r_2, r_1+r_2+r_3<\Delta} \!\!\!\!\!\!  dr_1 dr_2 dr_3
\frac{p_{\rho r_2}(0)p_{(1+\rho)\Delta-\rho(r_1+r_2+r_3)}(0)}{p_{(1+\rho)\Delta}(0)^2 (\log (r_2))^{2\xi}}  \\[-2ex]
&& \qquad \qquad \qquad \qquad \qquad  \times \quad \big|p_{\rho(r_1+r_3)}(0)p_{(1+\rho)\Delta}(0)-p_{\rho(r_1+r_2+r_3)}(0)p_{(1+\rho)\Delta-\rho r_2}(0)\big|  \\[2ex]
&\leq&
\frac{C\Delta^{\frac{5}{2}}}{\rho^\frac{3}{2}}   \iint\limits_{0<r_1, r_3<\Delta \atop A_1<r_2<\Delta} \!\!\! dr_1 dr_2 dr_3
\Big(\frac{p_{\rho(r_1+r_3)}(0)|p_{(1+\rho)\Delta}(0)-p_{(1+\rho)\Delta-\rho r_2}(0)|}{(1+r_2)^{\frac{3}{2}}(\log (r_2))^{2\xi}} \\
&& \qquad \qquad \qquad  \qquad \qquad \qquad \qquad + \quad
\frac{p_{(1+\rho)\Delta-\rho r_2}(0)|p_{\rho(r_1+r_3)}(0)-p_{\rho(r_1+r_2+r_3)}(0)|}{(1+r_2)^{\frac{3}{2}}(\log (r_2))^{2\xi}}\Big),
\end{eqnarray*}
where we applied (\ref{pDeltacomp}) to $\frac{p_{(1+\rho)\Delta-\rho(r_1+r_2+r_3)}(0)}{p_{(1+\rho)\Delta}(0)}$. Using Lemma \ref{L:est2}, we have
\begin{eqnarray*}
&& \frac{C \Delta^{\frac{5}{2}}}{\rho^{\frac{3}{2}}}  \iint\limits_{0<r_1, r_3<\Delta \atop A_1<r_2<\Delta} \frac{p_{\rho(r_1+r_3)}(0)|p_{(1+\rho)\Delta}(0)-p_{(1+\rho)\Delta-\rho r_2}(0)|}{(1+r_2)^{\frac{3}{2}}(\log (r_2))^{2\xi}}dr_1 dr_2 dr_3 \\
&\leq&
\frac{C \Delta^{\frac{5}{2}}}{\rho^3}  \iint\limits_{0<r_1, r_2, r_3<\Delta}  \frac{\frac{\rho r_2}{\Delta^{\frac{5}{2}}}}{(1+r_1+r_3)^{\frac{3}{2}}(1+r_2)^{\frac{3}{2}}(\log (e+r_2))^{2\xi}}dr_1 dr_2 dr_3 \\
&\leq& \frac{C}{\rho^2} \int_0^\Delta \frac{1}{\sqrt{r_2}}dr_2 \iint_{0\leq r_1+r_3\leq 2\Delta} \frac{1}{(r_1+r_3)^{\frac{3}{2}}}dr_1 dr_3 \\
&\leq& \frac{C\Delta}{\rho^2}.
\end{eqnarray*}
Similarly,
\begin{eqnarray*}
&& \frac{C\Delta^{\frac{5}{2}}}{\rho^{\frac{3}{2}}} \iint\limits_{0<r_1, r_3<\Delta \atop A_1<r_2<\Delta} \frac{p_{(1+\rho)\Delta-\rho r_2}(0)|p_{\rho(r_1+r_3)}(0)-p_{\rho(r_1+r_2+r_3)}(0)|}{(1+r_2)^{\frac{3}{2}}(\log (r_2))^{2\xi}}dr_1 dr_2 dr_3 \\
&\leq& \frac{C\Delta}{\rho^\frac{3}{2}}   \iint\limits_{0<r_1, r_3<\infty \atop A_1<r_2<\infty} \frac{\frac{ r_2}{\rho^{\frac{3}{2}}(r_1+r_3)^{\frac{3}{2}}(r_1+r_2+r_3)}}{(1+r_2)^{\frac{3}{2}}(\log (r_2))^{2\xi}}dr_1 dr_2 dr_3 \\
&\leq&  \frac{C \Delta}{\rho^3} \int_{A_1}^\infty \int_0^\infty \frac{1}{\sqrt{r_2}(\log(r_2))^{2\xi}}  \frac{1}{\sqrt{w}(r_2+w)}dw dr_2 \\
&=&  \frac{C \Delta}{\rho^3}\int_{A_1}^\infty \int_0^\infty \frac{1}{r_2(\log(r_2))^{2\xi}}  \frac{1}{\sqrt{t}(1+t)}dt dr_2 \\
&\leq& \frac{C \Delta}{\rho^3},
\end{eqnarray*}
where we used the fact that $\iint_{[0,\infty)^2} f(r_1+r_3)dr_1 dr_3 = \int_0^\infty wf(w)dw$, and made a change of variable $w=r_2 t$. Thus we have proved
$(3) \leq C\rho^{-3}\Delta$, which concludes the proof of (\ref{Hsigmavar}).
\qed
\bigskip

\appendix

\section{Renewal estimates}\label{S:renewal}

Consider a renewal process $\sigma:=\{\sigma_0=0<\sigma_1<\cdots\}$ on $[0,\infty)$, where
$(\sigma_i-\sigma_{i-1})_{i\in\N}$ are i.i.d.\ with distribution
$K(t)dt$ for a bounded density $K$ on $\R_+$ satisfying
\begin{equation}
\label{eq:Ktail}
K(t) \sim c_K t^{-1-\alpha}, \quad t\to\infty
\end{equation}
for some $\alpha \in (0,1)$ and $c_K \in (0,\infty)$. Let $K^{*n}$ denote the $n$-fold
convolution of $K$ with itself, and let $P(t) = \sum_{n=1}^\infty K^{*n}(t)$, as
defined in (\ref{Pt1}), be the corresponding renewal density.

We prove in Lemma \ref{L:renewalprob} a special case of the continuous time version of Doney's local limit theorem for renewal processes with infinite mean
\cite[Thm.~3]{D97}. Note that \cite{D97} allows a general regularly varying function in the right-hand side of (\ref{eq:Ktail}).
We stick to the narrower class, which suffices for our purposes, for the sake of a less cumbersome proof.

\bl\label{L:renewalprob} We have
\be\label{Psint}
\lim_{t\to\infty} c_K t^{1-\alpha} P(t) =
\frac{\alpha \sin(\alpha \pi)}{\pi}.
\ee \el

\bl\label{L:renewalcounts}
There exists a positive stable random variable 
$G$ with exponent $\alpha$, such that
\be
\lim_{t\to\infty} \P\big(|\sigma \cap [0,t]|\geq a t^\alpha\big) =
 \P\Big(G\leq \frac{1}{a^{1/\alpha}}\Big)
\qquad \mbox{for all $a> 0$}.
\ee
\el
It is well known that
\begin{equation}
\label{eq:stablelimit}
\frac{\sigma_n}{n^{\frac{1}{\alpha}}}=\frac{(\sigma_1-\sigma_0)+\cdots+(\sigma_n-\sigma_{n-1})}{n^{\frac{1}{\alpha}}}
\mathop{\to}^d G \quad \mbox{as} \; n\to\infty,
\end{equation}
where $G$ is a one-sided stable random variable of index $\alpha$.
Note that $\int_t^\infty K(s) \, ds \sim (c_K/\alpha) t^{-\alpha}$, thus
the normalisation is chosen here in such a way that
$\E[e^{-\lambda G}] = \exp\big(- \frac{c_K \Gamma(1-\alpha)}{\alpha}
\lambda^\alpha\big)$, $\lambda \geq 0$, i.e., $G$ is
$(c_K \Gamma(1-\alpha)/\alpha)^{1/\alpha}$ times a ``standard'' one-sided
stable random variable of index $\alpha$ (see, e.g.,
\cite[Thm.~XIII.6.2]{F66}). Since the characteristic function
of $G$ decays faster than any polynomial at infinity, $G$ has a $C^\infty$
density $g$, see, e.g., \cite[p.~48]{IL71}. As $G$ is a limit of
non-negative random variables, we must have $g(x)=0$ for $x<0$,
implying $g(0)=0$ by continuity. Furthermore, $g(x) \sim c_G x^{-1-\alpha}$
for $x\to\infty$ with some $c_G \in (0,\infty)$, see, e.g.,
\cite[Thm.~2.4.1]{IL71}.
In particular, $x \mapsto x^{-\alpha} g(x)$ is bounded and uniformly
continuous with $x^{-\alpha} g(x) \leq c/(1+x^{1+2\alpha})$
for some $c<\infty$. We have
\begin{equation}
\label{eq:minusalphamomentofG}
\int_0^{\infty} x^{-\alpha} g(x) \, dx
= \frac{\alpha}{c_K \Gamma(1-\alpha)} \cdot \frac{1}{\Gamma(1+\alpha)}
= \frac{1}{c_K \Gamma(1-\alpha) \Gamma(\alpha)}
= \frac{\sin(\alpha \pi)}{c_K \pi} .
\end{equation}
For the first equality, note that
$G^{-\alpha} = \Gamma(\alpha)^{-1} \int_0^\infty \lambda^{\alpha-1}
e^{-\lambda G} \, d\lambda$. The second identity uses well-known
facts about the $\Gamma$ function.
\smallskip

\noindent
{\bf Proof of Lemma~\ref{L:renewalcounts}.}
By (\ref{eq:stablelimit}),
$$
\lim_{t\to\infty} \P\big(|\sigma \cap [0,t]|\geq a t^\alpha\big)
= \lim_{t\to\infty} \P\Big(\frac{\sigma_{\lfloor a t^\alpha \rfloor}}{\lfloor a t^\alpha \rfloor^{1/\alpha}} \leq \frac{t}{\lfloor a t^\alpha \rfloor^{1/\alpha}}\Big)  = \P(G\leq \frac{1}{a^{1/\alpha}}),
$$
since the distribution of $G$ contains no atoms.
\qed

We will need the following uniform one-sided large
deviation estimate.
\begin{lem}
\label{lem:doneycts2}
We have for any sequence $c_n \to \infty$
\begin{equation}
\label{eq:doneycts2}
\lim_{n\to\infty} \sup_{t \geq c_n n^{\frac{1}{\alpha}}}
\Big| \frac{K^{* n}(t)}{n K(t)} - 1 \Big| = 0.
\end{equation}
\end{lem}
{\bf Proof.}
This follows from \cite[Thm.~1]{Z99} by specialising to the
one-dimensional asymmetric case. Note that Zaigraev\ \cite{Z99}
attributes the result in the present case (one-dimensional situation,
$K$ in the normal domain of attraction of a stable law) to
Tka\v{c}uk\ \cite{T73}, which the authors unfortunately could not
access.
\qed

\noindent
{\bf Proof of Lemma~\ref{L:renewalprob}.} Our proof follows more or less the scheme of
\cite[Thm.~3]{D97}, with \cite[Thm.~2]{D97} replaced by
Lemma~\ref{lem:doneycts2}. Even though we use Lemma~\ref{L:renewalprob}
in this paper only for $\alpha=1/2$, the proof is the same for all $\alpha \in (0,1)$.
\smallskip

By a local limit theorem for sums of random variables in the
domain of attraction of a stable law, e.g.\ \cite[Thm.~4.3.1]{IL71},
we have
\begin{equation}
\sup_{t \in \R_+} \big| n^{\frac{1}{\alpha}} K^{*n}(n^{\frac{1}{\alpha}} t) - g(t) \big|
\rightarrow 0 \quad \mbox{as}\;\; n\to\infty,
\end{equation}
where $g$ is the density of the one-sided stable random variable
appearing as the limit in (\ref{eq:stablelimit}).
Thus, we can find a continuous, strictly decreasing function
$\rho :[0,\infty) \to (0,\infty)$ with $\lim_{t\to\infty} \rho(t)=0$
such that
\begin{equation}
\sup_{t \in \R_+} \big| n^{\frac{1}{\alpha}} K^{*n}(n^{\frac{1}{\alpha}} t) - g(t) \big| \leq \rho(n)
\quad \mbox{for}\; n \in \N.
\end{equation}
Obviously, $\rho^{-1} : (0,\rho(0)] \to [0,\infty)$ is continuous and
strictly decreasing with $\lim_{y\to 0+} \rho^{-1}(y)=\infty$.
Note that the function $\psi : (0,\rho(0)^{1/(2-\alpha)}] \to [0,\infty)$ with
$\psi(y)=(\rho^{-1}(y^{2-\alpha}))^{1/\alpha}/y$ is strictly decreasing, and $\lim_{y\to 0+} \psi(y)= \infty$.
Define $\delta(t) := \psi^{-1}(t)$ for $t \geq 0$. Observe that then
$t \mapsto \delta(t)$ is strictly decreasing
and satisfies $\lim_{t\to\infty} \delta(t) = 0$.
Furthermore,
\begin{equation}
\rho\big((\delta(t) t)^\alpha\big) =
\rho\big(\big(\delta(t) \psi(\delta(t)) \big)^\alpha \big)
= \rho\big(\rho^{-1}\big(\delta(t)^{2-\alpha}\big)\big)
= \delta(t)^{2-\alpha},
\end{equation}
proving that $t \delta(t) \to \infty$ as $t\to\infty$, and
\begin{equation}
\label{eq:doneycontdeltaok}
\frac{\rho\big((\delta(t) t)^\alpha\big)}{\delta(t)^{1-\alpha}}
= \delta(t) \to 0 \quad \mbox{as} \; t \to \infty.
\end{equation}

Decompose
\begin{equation}
t^{1-\alpha} \sum_{n\geq 1} K^{*n}(t) =
t^{1-\alpha} \sum_{n > (\delta(t) t)^\alpha} K^{*n}(t) +
t^{1-\alpha} \sum_{n=1}^{[(\delta(t) t)^\alpha]} K^{*n}(t) =: S_1 + S_2.
\end{equation}
We have
\[
S_1 = t^{1-\alpha} \sum_{n > (\delta(t) t)^\alpha}
\frac{1}{n^{\frac{1}{\alpha}}} g\big(\frac{t}{n^{\frac{1}{\alpha}}}\big) +
t^{1-\alpha} \sum_{n > (\delta(t) t)^\alpha}
\frac{1}{n^{\frac{1}{\alpha}}} \Big( n^{\frac{1}{\alpha}} K^{*n}(n^{\frac{1}{\alpha}} \frac{t}{n^{\frac{1}{\alpha}}}\big) -
g\big(\frac{t}{n^{\frac{1}{\alpha}}}\big) \Big) =: S_1' +R_1,
\]
where
\[
|R_1| \leq \rho\big((\delta(t) t)^\alpha\big) \:
t^{1-\alpha} \hspace{-1em} \sum_{n> (\delta(t) t)^\alpha} \frac{1}{n^{\frac{1}{\alpha}}}.
\]
Since $\sum_{n> (\delta(t) t)^\alpha} \frac{1}{n^{\frac{1}{\alpha}}}
\sim \int_{(\delta(t) t)^\alpha}^\infty x^{-1/\alpha}\,dx
\sim \frac{\alpha}{1-\alpha} (\delta(t) t)^{\alpha-1}$, we obtain from
(\ref{eq:doneycontdeltaok}) that $R_1\to 0$ as $t\to\infty$.
\smallskip

Put $x^{(t)}_n := t/n^{\frac{1}{\alpha}}$, then we have
\[
\frac{t}{n^{\frac{1}{\alpha}}} \sim \alpha n \big( x^{(t)}_n - x^{(t)}_{n+1} \big)
= \alpha (t/x^{(t)}_n)^{\alpha} \big( x^{(t)}_n - x^{(t)}_{n+1} \big)
\]
since ${\alpha n n^{\frac{1}{\alpha}} \big( x^{(t)}_n - x^{(t)}_{n+1} \big)}/{t}
= \alpha n \big( 1 - \frac{n^{\frac{1}{\alpha}}}{(n+1)^{\frac{1}{\alpha}}}\big)
\sim \alpha n \big( 1 - \big(1-\frac{1}{n+1}\big)^{1/\alpha}\big) \to 1$, and hence
\[
S_1' \sim \alpha \sum_{n> (\delta(t) t)^\alpha} \frac{(t/x^{(t)}_n)^\alpha}{t^\alpha}
\big( x^{(t)}_n - x^{(t)}_{n+1} \big) g(x^{(t)}_n)
\sim \alpha \sum_{n\,:\, n^{\frac{1}{\alpha}} > \delta(t) t}
\big( x^{(t)}_n - x^{(t)}_{n+1} \big) \big(x^{(t)}_n\big)^{-\alpha} g(x^{(t)}_n).
\]
The term on the right is an approximating Riemann sum
and $n^{\frac{1}{\alpha}} > \delta(t) t$ means $x^{(t)}_n < 1/\delta(t)$, which tends to $\infty$ as
$t\to\infty$. Thus, recalling (\ref{eq:minusalphamomentofG}) and the discussion above it,
we have
\[
S_1 \to \alpha \int_0^{\infty} x^{-\alpha} g(x) \, dx
= \frac{\alpha \sin(\alpha \pi)}{c_K \pi}
\quad \mbox{as}\;\; t \to \infty.
\]
\smallskip

To bound $S_2$, note that $n \leq (\delta(t) t)^\alpha$ implies
$t \geq n^{1/\alpha}/\delta(t) \geq n^{1/\alpha}$ for $t$ sufficiently
large. In particular, for such $t$ and $n$, $\delta(t) \leq
\delta(n^{1/\alpha})$, so $t \geq (\delta(n^{1/\alpha}))^{-1} n^{1/\alpha}$.
Applying Lemma~\ref{lem:doneycts2} with $c_n := 1/\delta(n^{1/\alpha})
\to \infty$, we see that there exists $n_0 \in \N$, $t_0 < \infty$ and
$C< \infty$ such that
\begin{equation}
\label{eq:Kconvtail}
K^{*n}(t) \leq C n K(t) \quad \mbox{ for all}\;\:
n \geq n_0, \ t \geq \frac{n^{1/\alpha}}{\delta(t)} \vee t_0.
\end{equation}
Note that
\begin{equation}
K^{*n}(t) \leq 2n c_K (t/n)^{-1-\alpha} \leq 4 n^{2+\alpha} K(t)
\quad \mbox{for $t$ sufficiently large},
\end{equation}
which follows from (\ref{eq:Ktail}) and the observation that
$K^{*n}(t)$ is bounded from above by
\begin{eqnarray*}
\lefteqn{\sum_{j=1}^n \hspace{1em}
\idotsint\limits_{\sigma_0=0<\sigma_1<\cdots <\sigma_m=t}
1_{\{\sigma_j - \sigma_{j-1} \geq t/n\}}
\prod_{i=1}^m K(\sigma_i-\sigma_{i-1}) \prod_{i=1}^{m-1} d\sigma_i} \\
& \leq &
2c_K \Big(\frac{t}{n}\Big)^{-1-\alpha}
\sum_{j=1}^n \hspace{1em}
\idotsint\limits_{\sigma_0=0<\sigma_1<\cdots <\sigma_m=t}
\prod_{i=1, \, i \neq j}^m K(\sigma_i-\sigma_{i-1}) \prod_{i=1}^{m-1} d\sigma_i
\\
& \leq & 2 n c_K \Big(\frac{t}{n}\Big)^{-1-\alpha}
\bigg( \int_0^\infty K(\sigma)\,d\sigma \bigg)^{n-1}.
\end{eqnarray*}
Therefore if the constant $C$ appearing in (\ref{eq:Kconvtail})
is suitably increased, (\ref{eq:Kconvtail}) holds for all $n \in \N$.
Thus for $t$ sufficiently large, we have
\begin{eqnarray*}
S_2 & = & t^{1-\alpha} \sum_{n=1}^{[(\delta(t) t)^\alpha]} K^{*n}(t)
\leq 2 C c_K t^{-2\alpha} \sum_{n=1}^{[(\delta(t) t)^\alpha]} n
\leq 2 C c_K t^{-2\alpha} (\delta(t) t)^{2\alpha} = 2C c_K \delta(t)^{2\alpha},
\end{eqnarray*}
which converges to $0$ as $t\to\infty$.
\qed

\section{Random walk estimates}\label{S:RW}

\bl\label{L:LCLT}{\bf [Local central limit theorem]}
Let $(X_t)_{t\geq 0}$ with $X_0=0$ be a continuous time random walk on $\Z^d$ with jump rate $1$ and jump probability kernel $(q(x))_{x\in\Z^d}$, which is
irreducible and symmetric with finite covariance matrix $Q_{ij}=\sum_{x\in\Z^d} x_ix_j q(x)$, $1\leq i,j\leq d$. Let $p_t(\cdot)$ denote the transition probability kernel of $X$ at time $t$. Then
\be\label{LCLT0}
p_t(x) \leq p_t(0) \qquad \mbox{for all } x\in \Z^d \mbox{ and } t\geq 0,
\ee
and
\be\label{LCLT}
\lim_{t\to\infty} (2\pi t)^{\frac{d}{2}}\sqrt{{\rm det}\,Q}\, p_t(0) = 1.
\ee
\el
{\bf Proof.} Since $\hat p(k) :=\sum_{x\in\Z^d} e^{i\langle k, x\rangle} p_t(x) = e^{-t(1-\phi(k))}$, where $\phi(k)=\sum_{x\in\Z^d} e^{i\langle k,x\rangle}q(x)$
is real by the symmetry of $q$, by inverse Fourier transform,
$$
p_t(x) = \frac{1}{(2\pi)^d} \int_{[-\pi,\pi]^d} e^{-i\langle k, x\rangle} e^{-t(1-\phi(k))} dk \leq \frac{1}{(2\pi)^d} \int_{[-\pi,\pi]^d} e^{-t(1-\phi(k))} dk =p_t(0).
$$
For (\ref{LCLT}), see e.g.\ \cite[Prop.~7.9, Chap.~II]{S76} where a discrete time version was proved. The proof for the continuous
time version is identical.
\qed

\bl\label{L:est1} Let $X$, $q(\cdot)$, and $p_t(\cdot)$ be as in Lemma \ref{L:LCLT} without the symmetry assumption on $q$.
Then for any $a,b,c>0$, there exists some $C>0$ depending only on $q$ such that
\be
\sum_{x\in \Z^d} p_a(x) p_b(x) p_c(x) \leq \frac{C}{(1+ab+bc+ca)^{\frac{d}{2}}}.
\ee
\el
{\bf Proof.} Without loss of generality, assume that $a\geq b\geq c$. By the local central limit theorem, there exists $C_1>0$ such that uniformly in $t>0$ and $x\in \Z^d$, we have $ p_t(x) \leq \frac{C_1}{(1+t)^{\frac{d}{2}}}$. Then
$$
\sum_{x\in\Z^d} p_a(x) p_b(x) p_c(x) \leq \frac{C_1^2}{(1+ab)^{\frac{d}{2}}} \sum_{x\in\Z^d} p_c(x) = \frac{C_1^2}{(1+ab)^{\frac{d}{2}}} \leq \frac{C}{(1+ab+bc+ca)^{\frac{d}{2}}}.
$$
\qed

\bl\label{L:est2} Let $X$, $q(\cdot)$, $Q$, and $p_t(\cdot)$ be as in Lemma \ref{L:LCLT} so that $q$ is symmetric.
Then there exist $C_1, C_2>0$ depending on $q$, such that
\be\label{est2eq1}
\frac{C_1 r}{t^{\frac{d}{2}}(t+r)} \leq p_t(0)-p_{t+r}(0) \leq \frac{C_2 r}{t^{\frac{d}{2}}(t+r)},
\ee
where the first inequality holds for all $r>0$, $t>1$, and the second inequality holds for all $r, t>0$.
\el
{\bf Proof.} By the symmetry of $q$, $\phi(k) := \sum_{x} e^{i\langle k,x\rangle} q(x) \in [-1,1]$, and $\E[e^{i\langle k,X_t\rangle}] = e^{-t(1-\phi(k))}$.
Therefore,
$$
p_t(0)-p_{t+r}(0) = \frac{1}{(2\pi)^d} \int_{[-\pi,\pi]^d} \big(e^{-t(1-\phi(k))}-e^{-(t+r)(1-\phi(k))}\big) dk.
$$
By irreducibility of $q(\cdot)$, $\phi(k)=1$ only at $k=0$, and hence $c:=\inf_{|k|\geq \eps, k\in[-\pi,\pi]^d} (1-\phi(k))>0$ for any $\eps>0$.
By Taylor expansion, if $\eps>0$ is sufficiently small, then
$$
\frac{1}{4}\langle k, Qk\rangle \leq (1-\phi(k)) \leq \langle k, Qk\rangle  \qquad \forall \, |k|<\eps.
$$
Therefore
\begin{eqnarray}
&& (2\pi)^d(p_t(0)-p_{t+r}(0)) = \int_{[-\pi,\pi]^d} e^{-t(1-\phi(k))}(1-e^{-r(1-\phi(k))}) dk  \nn \\
&\leq& r \int_{[-\pi,\pi]^d} (1-\phi(k)) e^{-t(1-\phi(k))} dk  \nn \\
&\leq& 2 r \int_{|k|> \eps, k\in [-\pi,\pi]^d} e^{-t(1-\phi(k))} dk + r \int_{|k|\leq \eps} \langle k, Qk\rangle e^{-\frac{t \langle k, Qk\rangle}{4}} dk  \nn \\
&\leq& 2 (2\pi)^d e^{-ct} r + \frac{r}{t^{\frac{d}{2}+1}} \int_{\R^d} \langle k, Qk\rangle e^{-\frac{\langle k, Qk\rangle}{4}} dk \nn \\
&\leq& \frac{C r}{t^{\frac{d}{2}+1}}, \label{lem3est}
\end{eqnarray}
which implies that $p_t(0)-p_{t+r}(0) \leq \frac{C_2 r}{t^{\frac{d}{2}}(t+r)}$ for $r<t$. When $r\geq t$, the same bound follows from the local central limit theorem.

Similarly,
\begin{eqnarray*}
&& (2\pi)^d(p_t(0)-p_{t+r}(0)) = \int_{[-\pi,\pi]^d} e^{-(t+r)(1-\phi(k))}(e^{r(1-\phi(k))}-1) dk \\
&\geq& r \int_{|k|\leq \eps, k\in [-\pi,\pi]^d} (1-\phi(k)) e^{-(t+r)(1-\phi(k))} dk \geq
r \int_{|k|\leq \eps, k\in [-\pi,\pi]^d} \frac{\langle k, Qk\rangle}{4} e^{-(t+r)\langle k, Qk\rangle}  dk \\
&\geq& \frac{Cr}{(t+r)^{\frac{d}{2}+1}},
\end{eqnarray*}
which follows by a change of variable for $k$ and the fact that $t+r>1$. This implies $p_t(0)-p_{t+r}(0) \geq \frac{C_1 r}{t^{\frac{d}{2}}(t+r)}$ for $r<t$.
When $r>t$, the same bound follows from the local central limit theorem.
\qed

\bl\label{L:est3} Let $X$, $q(\cdot)$ and $p_t(\cdot)$ be as in Lemma \ref{L:LCLT} so that $q$ is symmetric. Then there exist $C>0$ depending only on $q$ such that, for all $a,b>0$ and $t>0$,
\be\label{est3eq1}
|p_t(0)p_{t+a+b}(0)-p_{t+a}(0)p_{t+b}(0)| \leq \frac{C ab}{t^d(t+a)(t+b)}.
\ee
\el
{\bf Proof.} Note that
\begin{eqnarray}
\!\!\!\!\!\!\!\!\!\!\!\!\!\!\!\!\!\!\!\!&& p_t(0)p_{t+a+b}(0)-p_{t+a}(0)p_{t+b}(0)  \nonumber \\
\!\!\!\!\!\!\!\!\!\!\!\!\!\!\!\!&=&\!\!\!\! p_{t+a+b}(0)(p_t(0)-p_{t+a}(0)) - p_{t+a}(0)(p_{t+b}(0)-p_{t+a+b}(0))  \nonumber \\
\!\!\!\!\!\!\!\!\!\!\!\!\!\!\!\!&=&\!\!\!\! \big(p_{t+a+b}(0)-p_{t+a}(0)\big)\big(p_t(0)-p_{t+a}(0)\big) + p_{t+a}(0)( p_t(0)-p_{t+a}(0)-p_{t+b}(0)+p_{t+a+b}(0)).
\label{est3eq2}
\end{eqnarray}
By Lemma \ref{L:est2}, the first term in (\ref{est3eq2}) is bounded in absolute value by
$$
\frac{Cb}{(t+a)^{\frac{d}{2}}(t+a+b)}\cdot \frac{Ca}{t^{\frac{d}{2}}(t+a)},
$$
which is clearly bounded by the RHS of (\ref{est3eq1}).

For the second term in (\ref{est3eq2}), we claim that
\be\label{est3eq3}
0\leq p_t(0)-p_{t+a}(0)-p_{t+b}(0)+p_{t+a+b}(0) \leq \frac{Cab}{t^{\frac{d}{2}}(t+a)(t+b)},
\ee
which together with the fact that $p_{t+a}(0) \leq C t^{-\frac{d}{2}}$ imply (\ref{est3eq1}). Note that
\begin{eqnarray*}
(2\pi)^d \big(p_t(0)-p_{t+a}(0)-p_{t+b}(0)+p_{t+a+b}(0)\big) \! &=&  \!\!\!\! \int_{[-\pi,\pi]^d} \!\!\!\!\!\!\!\!\!\!\!
e^{-t(1-\phi(k))}(1-e^{-a(1-\phi(k))})(1-e^{-b(1-\phi(k))})dk \\
&\leq& \! ab \int_{[-\pi,\pi]^d} \!\!\!\! (1-\phi(k))^2 e^{-t(1-\phi(k))} dk.
\end{eqnarray*}
Clearly $p_t(0)-p_{t+a}(0)-p_{t+b}(0)+p_{t+a+b}(0)\geq 0$. For the upper bound, exactly as in (\ref{lem3est}), we can Taylor expand $\phi(k)$ around $k=0$ for $|k|\leq \eps$ and bound $|\phi(k)|$ uniformly
for $|k|>\eps$, which gives
$$
p_t(0)-p_{t+a}(0)-p_{t+b}(0)+p_{t+a+b}(0) \leq  \frac{Cab}{t^{\frac{d}{2}+2}}.
$$
When $a,b<t$, this implies (\ref{est3eq3}). If $b>t$, then (\ref{est3eq3}) follows from the bound
$$
p_t(0)-p_{t+a}(0)-p_{t+b}(0)+p_{t+a+b}(0) \leq \frac{Ca}{t^{\frac{d}{2}}(t+a)} + \frac{Ca}{(t+b)^{\frac{d}{2}}(t+a+b)}
$$
by Lemma \ref{L:est2}. The same argument applies when $a>t$.
\qed

\bl\label{L:retprocomp}{\bf [Comparison of return probabilities]}
Let $X$, $q(\cdot)$ and $p_t(\cdot)$ be as in Lemma \ref{L:LCLT} so that $q$ is symmetric. For $1\leq i\leq n$, let $a_i, b_i> 0$, and let $Z_i$ be an independent random variable distributed as $X_{a_i}$ conditioned on $X_{a_i+b_i}=0$. Then
\be\label{rtcomp}
\P(Z_1+\cdots+ Z_n =0) > \P(X_{a_1+\cdots + a_n}=0).
\ee
\el
{\bf Proof.} Let $\phi(k)=\sum_{x}e^{i\langle k,x\rangle}q(x)$ and $\psi_i(k) = \E[e^{i\langle k, Z_i\rangle}]$. Since
$\E[e^{i\langle k, X_t\rangle}] = e^{-t(1-\phi(k))}$, by Fourier transform, (\ref{rtcomp}) is equivalent to
\be\label{rtcomp2}
\int_{[-\pi,\pi]^d} \psi_1(k)\cdots \psi_n(k) {\rm d}k > \int_{[-\pi,\pi]^d} e^{-\sum_{i=1}^n a_i (1-\phi(k))} {\rm d}k.
\ee
By symmetry of $q$, $\phi(k) \in [-1,1]$ and $e^{- a_i(1-\phi(k))} \in (0,1]$. Therefore to verify (\ref{rtcomp2}), it suffices to show that
for each $1\leq i\leq n$,
\be\label{rtcomp2.5}
\psi_i(k) \geq e^{- a_i(1-\phi(k))}
\ee
for all $k\in [-\pi,\pi]^d$, with strict inequality for some $k\in [-\pi,\pi]^d$ .

Note that $\hat p_s(k):=\sum_x e^{i\langle k,x\rangle}p_s(x)=e^{-s(1-\phi(k))}$. By definition, $\P(Z_i=x) = \frac{p_{a_i}(x)p_{b_i}(x)}{p_{a_i+b_i}(0)}$, and hence
$$
\psi_i(k) = \frac{(\hat p_{a_i}*\hat p_{b_i})(k)}{p_{a_i+b_i}(0)} = \frac{\int_{[-\pi,\pi]^d} e^{-a_i(1-\phi(k-u))-b_i(1-\phi(u))}{\rm d}u}{\int_{[-\pi,\pi]^d} e^{-(a_i+b_i)(1-\phi(u))}{\rm d}u}.
$$
By symmetry, $\psi_i(k)=\psi_i(-k)$, and hence
\begin{eqnarray}
\psi_i(k) &=& \frac{\int_{[-\pi,\pi]^d} e^{-b_i(1-\phi(u))} \frac{e^{-a_i(1-\phi(k-u))}+e^{-a_i(1-\phi(-k-u))}}{2}{\rm d}u}{\int_{[-\pi,\pi]^d} e^{- (a_i+b_i)(1-\phi(u))}{\rm d}u}  \nonumber \\
&\geq& \frac{\int_{[-\pi,\pi]^d} e^{- b_i(1-\phi(u))} e^{- a_i(1-\frac{\phi(k-u)+\phi(-k-u)}{2})}{\rm d}u}{\int_{[-\pi,\pi]^d} e^{-(a_i+b_i)(1-\phi(u))}{\rm d}u},
\label{rtcomp3}
\end{eqnarray}
where we applied Jensen's inequality. Note that since $\phi(x)$ is not identically equal to $1$, for some choice of $k$ and $u$, we have $\phi(k-u)\neq\phi(-k-u)$ so that there is strict inequality in (\ref{rtcomp3}) for some $k$. By symmetry,
\begin{eqnarray}
\phi(k-u)+ \phi(-k-u) &=& \sum_x q(x) (e^{i\langle k-u,x\rangle}+ e^{i\langle -k-u,x\rangle}) \nonumber \\
&=& \sum_x q(x)\big(\cos \langle k-u,x\rangle + \cos \langle -k-u,x\rangle\big) \nonumber \\
&=& 2\sum_x q(x)\cos \langle k,x\rangle \cos\langle u,x\rangle  \nonumber \\
&\geq& 2 \sum_x q(x) \big(\cos\langle k,x\rangle + \cos\langle u,x\rangle -1 \big) \nonumber\\
&=& 2(\phi(k)+\phi(u)-1),  \label{rtcomp4}
\end{eqnarray}
where we used $(1-\cos\alpha)(1-\cos\beta)\geq 0$. Plugging this bound into (\ref{rtcomp3}) then yields (\ref{rtcomp2.5}).
\qed

\section{Proof of Theorem~\ref{T:monshif}}\label{S:monshif}
Let $\rho'>\rho\geq 0$. Let $X, Y, Y^{(1)}, Y^{(2)}$ be independent random walks on $\Z^d$ with the same symmetric jump kernel with finite second moments and with respective jump rates $1, \rho, \frac{1+\rho'}{1+\rho} \rho$ and $\frac{\rho'-\rho}{1+\rho}$. Then $Y':=Y^{(1)}+Y^{(2)}$ and $X':=X-Y^{(2)}$ are random walks
with the same jump kernel and jump rates $\rho'$ and $\frac{1+\rho'}{1+\rho}$, where for $X'$ we used the symmetry of the kernel. The key observation is that
\be\label{Ceq1}
\big(\E^{Y^{(2)}}[Z^\beta_{t,Y'}]\big)_{t>0} \stackrel{\rm law}{=} \Big(Z^{\beta\frac{1+\rho}{1+\rho'}}_{t\frac{1+\rho'}{1+\rho},Y}\Big)_{t>0},
\ee
which is a simple consequence of the fact that
\begin{eqnarray*}
\E^{Y^{(2)}}[Z^\beta_{t,Y'}] &=& \E^{Y^{(2)}, X}\big[e^{\beta \int_0^t 1_{\{X_s=Y^{(1)}_s+Y^{(2)}_s\}}{\rm d}s}\big] = \E^{X'}\big[e^{\beta \int_0^t 1_{\{X'_s=Y^{(1)}_s\}}{\rm d}s}\big], \\
Z^{\beta\frac{1+\rho}{1+\rho'}}_{t\frac{1+\rho'}{1+\rho},Y} &=& \E^X\Big[e^{\beta\frac{1+\rho}{1+\rho'} \int_0^{t\frac{1+\rho'}{1+\rho}} 1_{\{X_s=Y_s\}}{\rm d}s} \Big],
\end{eqnarray*}
and the fact that
$$
(X_{\frac{1+\rho'}{1+\rho}s}, Y_{\frac{1+\rho'}{1+\rho}s})_{s\geq 0} \stackrel{\rm law}{=} (X'_s, Y^{(1)}_s)_{s\geq 0}.
$$

Note that
$$
\lim_{t\to\infty} \frac{1}{t} \log \E^{Y^{(2)}}[Z^\beta_{t,Y'}] \geq \lim_{t\to\infty} \frac{1}{t} \log Z^\beta_{t,Y'} = F(\beta, \rho') \qquad a.s.
$$
On the other hand, by (\ref{Ceq1}),
$$
\lim_{t\to\infty} \frac{1}{t} \log \E^{Y^{(2)}}[Z^\beta_{t,Y'}] = \lim_{t\to\infty} \frac{1}{t} \log Z^{\beta\frac{1+\rho}{1+\rho'}}_{t\frac{1+\rho'}{1+\rho},Y}
= \frac{1+\rho'}{1+\rho} F\big(\beta\frac{1+\rho}{1+\rho'}, \rho\big) \qquad a.s.
$$
Therefore
\be
F(\beta, \rho') \leq \frac{1+\rho'}{1+\rho} F\big(\beta\frac{1+\rho}{1+\rho'}, \rho\big) \qquad \mbox{for all } \rho'>\rho\geq 0,
\ee
which implies the first inequality in (\ref{monoton1}).

Similarly, by (\ref{Ceq1}),
$$
\sup_{t>0} Z^{\beta\frac{1+\rho}{1+\rho'}}_{t\frac{1+\rho'}{1+\rho},Y} <\infty \quad a.s. \quad \Longleftrightarrow \quad
\sup_{t>0} \E^{Y^{(2)}}[Z^\beta_{t,Y'}] <\infty \quad a.s.  \quad \Longrightarrow \quad \sup_{t>0} Z^\beta_{t,Y'}<\infty \quad a.s.,
$$
which implies the second inequality in (\ref{monoton1}).

To prove the first inequality in (\ref{monoton2}), let $\phi(\rho):=\frac{\beta_c(\rho)}{1+\rho}$, recall that $\beta_c^{\rm ann}(\rho) = (1+\rho)/G$, and note that
\begin{eqnarray*}
\beta_c(\rho')-\beta^{\rm ann}_c(\rho') - (\beta_c(\rho)-\beta^{\rm ann}_c(\rho)) &=& (1+\rho')(\phi(\rho')-G^{-1}) - (1+\rho)(\phi(\rho)-G^{-1}) \\
&=& (\rho'-\rho)(\phi(\rho')-G^{-1}) + (1+\rho)(\phi(\rho')-\phi(\rho))>0,
\end{eqnarray*}
since $\phi(\rho)$ is non-decreasing in $\rho$, and $(1+\rho')(\phi(\rho')-G^{-1})=\beta_c(\rho')-\beta^{\rm ann}_c(\rho')>0$ by Theorem~\ref{T:cptcts}
and its analogue in dimensions $d\geq 4$ shown in \cite{BS09}. The proof of the second inequality in (\ref{monoton2}) is identical.
\qed
\bigskip

\noindent
{\bf Acknowledgment} We thank F.L.~Toninelli for sending us the preprint \cite{BT09} before publication. We thank the referee for a careful reading of the paper and
helpful comments, and in particular, for bringing to our attention Theorem~\ref{T:monshif} and allowing us to include its elegant proof here.


\begin{thebibliography}{}


\bibitem[BT09]{BT09}
Q.\ Berger, F.L.\ Toninelli.
On the critical point of the Random Walk Pinning Model in dimension $d=3$,
arXiv:0911.1661v2, 2009. To appear in {\em Electron.\ Journal Probab.}


\bibitem[BGdH08]{BGdH08}
M.\ Birkner, A.\ Greven, F.\ den Hollander.
Collision local time of transient random walks and intermediate phases in
interacting stochastic systems,
EURANDOM Report 2008-49, 2008.
\url{http://www.eurandom.nl/reports/2008/049-report.pdf}

\bibitem[BS09]{BS09}
M.\ Birkner, R.\ Sun.
Annealed vs quenched critical points for a random walk pinning model, arXiv:0807.2752v2, 2009.
To appear in {\em Ann.\ Inst.\ Henri Poincar\'e  Probab.\ Stat}.

\bibitem[DGLT09]{DGLT09}
B.\ Derrida, G.\ Giacomin, H.\ Lacoin, F.L.\ Toninelli.
Fractional moment bounds and disorder relevance for pinning models,
{\em Commun.\ Math.\ Phys.} 287, 867--887, 2009.

\bibitem[D97]{D97}
R.A.\ Doney.
One-sided local large deviation and renewal theorems in the case of infinite mean,
{\em Probab.\ Theory Rel.\ Fields} 107, 451--465, 1997.

\bibitem[D96]{D96}
R.\ Durrett. {\it Probability: Theory and Examples}, 2nd ed., Duxbury Press, 1996.

\bibitem[F66]{F66}
W.\ Feller.
{\em An introduction to probability theory and its applications}. Vol.\ II.
John Wiley \& Sons, Inc., New York-London-Sydney, 1966.

\bibitem[G07]{G07}
G.\ Giacomin.
{\it Random Polymer Models}, Imperial College Press, World Scientific, 2007.

\bibitem[GLT08]{GLT08}
G.\ Giacomin, H.\ Lacoin and F.L.\ Toninelli.
Marginal relevance of disorder for pinning models,
{\it Commun.\ Pure Appl.\ Math} 63, 233--265, 2010.

\bibitem[GLT09]{GLT09}
G.\ Giacomin, H.\ Lacoin and F.L.\ Toninelli.
Disorder relevance at marginality and critical point shift, arXiv:0906.1942v1, 2009.
To appear in {\em Ann.\ Inst.\ Henri Poincar\'e  Probab.\ Stat}.

\bibitem[GdH07]{GdH07}
A.\ Greven and F.\ den Hollander. Phase transitions for the long-time behaviour of interacting diffusions,
{\it Ann.\ Probab.\ } 35, 1250--1306, 2007.

\bibitem[IL71]{IL71}
I.A.\ Ibragimov, Yu.V.\ Linnik,
\emph{Independent and stationary sequences of random variables}.
Wolters-Noordhoff Publishing, Groningen, 1971.

\bibitem[L09]{L09}
H.\ Lacoin.
New bounds for the free energy of directed polymers in dimension $1+1$ and $1+2$,
{\em Commun.\ Math.\ Phys.} 294, 471--503, 2010.

\bibitem[S76]{S76}
F.\ Spitzer.
{\em Principles of random walks}, 2nd ed. Springer-Verlag, New York-Heidelberg, 1976.

\bibitem[T73]{T73}
S.G.\ Tka\v{c}uk,
Local limit theorems, allowing for large deviations,
in the case of stable limit laws. (Russian.\ Uzbek summary)
\emph{Izv.\ Akad.\ Nauk UzSSR Ser.\ Fiz.-Mat.\ Nauk} 17,
no.\ 2, 30--33, 70, 1973.

\bibitem[YZ09]{YZ09}
A.~Yilmaz, O.~Zeitouni.
Differing averaged and quenched large deviations for random walks in random environments in dimensions two and three,
arXiv:0910.1169v1, 2009. To appear in {\em Commun.\ Math.\ Phys.}

\bibitem[Z99]{Z99}
A.\ Zaigraev,
Multivariate large deviations with stable limit laws,
\emph{Probab.\ Math.\ Statist.} 19,  no.\ 2,
Acta Univ.\ Wratislav.\ No.\ 2198, 323--335, 1999.


\end{thebibliography}
\end{document}